\documentclass[10pt]{iopart}

\usepackage{iopams}
\usepackage[colorlinks, linkcolor=blue,citecolor=red]{hyperref} 
\usepackage{enumerate}
\usepackage[T1]{fontenc}
\usepackage{aecompl}
\usepackage{geometry}
\usepackage[normalem]{ulem}
\useunder{\uline}{\ul}{}
\usepackage{mathrsfs}
\usepackage{epsfig}
\usepackage{color}
\usepackage{cite}
\usepackage{amsthm}
\usepackage{bm}
\usepackage{amssymb}
\usepackage{multirow}
\usepackage{mathrsfs}
\usepackage{array}
\usepackage{algorithm}
\usepackage{algorithmic}

\usepackage{caption}
\newtheorem{thm}{Theorem}[section]
\newtheorem{definition}[thm]{Definition}
\newtheorem{lemma}[thm]{Lemma}
\newtheorem{corollary}[thm]{Corollary}
\newtheorem{prop}[thm]{Proposition}
\newtheorem{assumption}{Assumption}[section]
\newtheorem{remark}{Remark}[section]
\newtheorem{alg}{Algorithm}[section]

\makeatletter 

\@addtoreset{equation}{section}
\makeatother
\usepackage{graphicx}
\usepackage{subfigure}
\usepackage{color}
\usepackage{multirow}
\usepackage{epstopdf}
\usepackage{lineno}

\usepackage{indentfirst}
\usepackage{footnote}
\usepackage{threeparttable}
\begin{document}
\captionsetup{labelformat=default,labelsep=period}
\title[]{A fast two-point gradient algorithm based on sequential subspace optimization method for nonlinear ill-posed problems}

\author{Guangyu Gao$^1$, Bo Han$^1$\footnote{E-mail: bohan@hit.edu.cn} and Shanshan Tong$^2$}


\vspace{10pt}

\address{$^{1}$ Department of Mathematics, Harbin Institute of Technology, Harbin, Heilongjiang Province, 150001, PR China}
\address{$^{2}$ School of Mathematics and Information Science, Shaanxi Normal University, Xi'an, Shaanxi Province, 710119, PR China}

\ead{guangyugao60@163.com, bohan@hit.edu.cn, tongshanshan33@163.com}


\begin{abstract}
 In this paper, we propose and analyze a fast two-point gradient algorithm for solving nonlinear ill-posed problems, which is based on the sequential subspace optimization method. A complete convergence analysis is provided under the classical assumptions for iterative regularization methods.
 The design of the two-point gradient method involves the choices of the combination parameters which is systematically discussed. Furthermore, detailed numerical simulations are presented for inverse potential problem, which exhibit that the proposed method leads to a strongly decrease of the iteration numbers and the overall computational time can be significantly reduced.

\vspace{2pc}
\end{abstract}

\section{Introduction}\label{section-1}
\noindent In this paper, we mainly focus on an iterative solution of nonlinear inverse problems in Hilbert spaces which can be written as
  \begin{equation}\label{equation 1-1}
    F(x)=y.
  \end{equation}
Here $F:\mathcal{D}(F)\subset \mathcal{X}\rightarrow \mathcal{Y}$, is possibly nonlinear operator between real Hilbert spaces $\mathcal{X}$ and
$\mathcal{Y}$, with domain of definition $\mathcal{D}(F)$. The solution set of this problem is defined as
\[
M_{F(x)=y}:=\left\{x\in\mathcal{X}:F(x)=y\right\}.
\]
Instead, only approximate measured data $y^\delta$ is available such that
  \begin{equation}\label{equation 1-2}
   \left\|y^\delta-y\right\|\leq\delta
  \end{equation}
  with $\delta>0$ (noise level).
  Due to the ill-posed (or ill-conditioned) nature of problem (\ref{equation 1-1}), regularization strategy is necessary to obtain their stable and accurate numerical solutions, and many effective techniques have been proposed over the past few decades (see e.g. \cite{1,2,3,4,5}).
  There are at least two common regularization methods for solving (\ref{equation 1-1}), i.e., (generalized) \emph{Tikhonov-type regularization} and  \emph{Iterative-type regularization}.

 In \emph{Tikhonov-type regularization}, one attempts to approximate an $x_0$-minimum-norm solution of  (\ref{equation 1-1}). The corresponding study of linear and nonlinear ill-posed problems is discussed in \cite{6}.  Furthermore, \emph{Iterative-type regularization} represents a very powerful and popular class of numerical solvers. Many publications have been concerned with iterative regularization methods for nonlinear problems in the Hilbert spaces.

 The most basic iterative method for solving (\ref{equation 1-1}) is the Landweber iteration \cite{1,7}, which reads
 \begin{equation}\label{equation 1-3}
 x_{k+1}^\delta=x_{k}^\delta-\alpha_{k}^\delta F'(x_{k}^\delta)^*\left(F(x_{k}^\delta)-y^\delta\right),~~k=0,1,\ldots,
  \end{equation}
  there $F'(x)^*$ is the Hilbert space adjoint of the derivative of $F$, $x_{0}^\delta=x_{0}$ is a suitable guess and $\alpha_{k}^\delta$ is a scaling parameter.  There are many definitions about $\alpha_{k}^\delta$ in the convergence analysis called the \emph{constant method}, \emph{steepest descent method} \cite{8} and the \emph{minimal error method} \cite{2}.

   It is well-known that Landweber iteration, the steepest descent method as well as minimal error method are quite slow. Hence, acceleration strategies have to be used in order to speed them up to make them more applicable in practise. Nesterov acceleration scheme was first  introduced for solving convex programming problems which is effective for implementation and greatly improves the convergence rate \cite{9,10,11}. Recently, Nesterov acceleration scheme was used to accelerate Landweber iteration in \cite{12,19}, leading to the following method
\begin{equation}\label{equation 1-4}
    \begin{array}{ll}
 z_{k}^\delta=x_{k}^\delta+\frac{k-1}{k+\alpha-1}(x_{k}^\delta-x_{k-1}^\delta), \\
 x_{k+1}^\delta=z_{k}^\delta-\alpha_{k}^\delta s_{k}^\delta,~~ s_{k}^\delta=F'(z_{k}^\delta)^*\left(F(z_{k}^\delta)-y^\delta\right),\\
 x_{0}^\delta=x_{-1}^\delta=x_{0},
    \end{array}
\end{equation}
where the constant $\alpha\geq3$. Although the convergence analysis of (\ref{equation 1-4}) is not given, the significant acceleration effect can be verified by numerical experiments \cite{12}. Furthermore, the two-point gradient (TPG) method was proposed in \cite{13}, and is given by
\begin{equation}\label{equation 1-5}
    \begin{array}{ll}
 z_{k}^\delta=x_{k}^\delta+\lambda_{k}^\delta(x_{k}^\delta-x_{k-1}^\delta), \\
 x_{k+1}^\delta=z_{k}^\delta-\alpha_{k}^\delta s_{k}^\delta,~~ s_{k}^\delta=F'(z_{k}^\delta)^*\left(F(z_{k}^\delta)-y^\delta\right),\\
 x_{0}^\delta=x_{-1}^\delta=x_{0}.
    \end{array}
\end{equation}
The convergence results are obtained under a suitable choice of the combination parameters $\{\lambda_{k}^\delta\}$. Recently, the study of TPG method for solving problem (\ref{equation 1-1}) in Banach spaces has been proposed in \cite{14}.

 Sequential subspace optimization (SESOP) method was first introduced to solve linear systems of equations in finite dimensional vectors spaces and then was brought up for dealing with nonlinear inverse problems in Hilbert spaces and Banach spaces \cite{15,16,17,18}. The significant improvement of SESOP method is to change one search direction to multiple as it is the case for Landweber or the conjugate gradient method. The general iteration format of SESOP is
 \begin{equation}\label{equation 1-6}
    \begin{array}{ll}
 x_{k+1}^\delta=x_{k}^\delta-\sum\limits_{i\in I_{k}^\delta}t_{k,i}^\delta F'(x_{i}^\delta)^*w_{k,i}^\delta, \;\;w_{k,i}^\delta=F(x_{i}^\delta)-y^\delta,\\
 x_{0}^\delta=x_{-1}^\delta=x_{0},
    \end{array}
\end{equation}
where $t_{k,i}^\delta$ is chosen such that $x_{k+1}^\delta$ intersects in the stripes or hyperplane of each $x_{i}^\delta$ which depend on search directions, forward operator and noise level, $I_{k}^\delta$ is the index set of the selected space.

 Motivated by the above considerations, an accelerated TPG method based on SESOP method (TGSS) is proposed for solving (\ref{equation 1-1}) in this paper, which leads to the following scheme
 \begin{equation}\label{equation 1-7}
    \begin{array}{ll}
 z_{k}^\delta=x_{k}^\delta+\lambda_{k}^\delta(x_{k}^\delta-x_{k-1}^\delta), \\
 x_{k+1}^\delta=z_{k}^\delta-\sum\limits_{i\in I_{k}^\delta}t_{k,i}^\delta F'(z_{i}^\delta)^*w_{k,i}^\delta, \quad w_{k,i}^\delta=F(z_{i}^\delta)-y^\delta,\\
 x_{0}^\delta=x_{-1}^\delta=z_{0}^\delta=x_{0},
    \end{array}
\end{equation}
where $t_{k,i}^\delta$ is chosen so that $x_{k+1}^\delta$ intersects in the stripes or hyperplane of each $z_{i}^\delta$. $I_{k}^\delta$ is the index set of the selected space that needs to contain $k$, i.e., the current search direction $F'(z_{k}^\delta)^*w_{k,k}^\delta$ need to be chosen
as one of the search directions, which guarantees the decent property of the iteration error. In order to find nontrivial $\lambda_{k}^\delta$, we apply the Nesterov acceleration scheme in \cite{9} and the discrete backtracking search (DBTS) algorithm in \cite{13} to our situation. The convergence analysis of the proposed method is provided with the general assumptions for iteration regularization methods. Furthermore, we provide a complete convergence analysis by showing a uniform convergence result for the noise-free case with the combination parameters chosen in a certain range \cite{14}.

 However, the TGSS method for solving (\ref{equation 1-1}), to the best of authors knowledge, have not been investigated before. On the basis of the iteration format of TGSS method, we summarize the following improvements. On the one hand, the method we present is to use finite search directions in each step of the iteration instead of just one as it is the case for the TPG method. This leads to more computational time in each step of the iteration, but it may reduce the total number of steps to obtain a satisfying approximation. On the other hand, compare to SESOP method, in each iteration, replacing $z_{i}^\delta$ with $x_{i}^\delta$ not only optimizes the search directions for each iteration but improves the initial projection point.

 An outline of the remaining paper is as follows. In section \ref{section-2}, we present the assumptions and required notations, and derive some preliminary results. In section \ref{section-3}, we introduce the TGSS method for nonlinear ill-posed problems with exact data and noisy data respectively, and give an algorithm with two search directions. The convergence and regularity analysis of TGSS are then given in section \ref{section-4}. In section \ref{section-5}, some numerical experiments are performed to present the efficiency of proposed methods. Finally, we give some conclusions in section \ref{section-6}.

\section{Assumptions and preliminary results}\label{section-2}
\noindent In this section, we state some main assumptions which are standard in the analysis of regularization methods \cite{20,21}. Furthermore, we give some basic definitions as well as their preliminary results, which we refer to \cite{15,16} for details.

\begin{assumption} \label{assumption 2-1}
Let $\rho$ be a positive number such that $B_{4\rho}(x_0)\subset\mathcal{D}(F)$, where $B_{4\rho}(x_0)$ denotes the closed ball around $x_0$ with radius $4\rho$.
The mapping $x\mapsto F'(x)$,  $F':\mathcal{D}(F)\rightarrow\mathcal{L}(\mathcal{X},\mathcal{Y})$ is a Fr\'{e}chet derivative operater.
\begin{description}
  \item[(1)] The equation $F(x)=y$ has a solution $x_*$ in $B_{\rho}(x_0)=B_{\rho}(x_{-1})$ which is not necessarily unique.
  \item[(2)] The local tangential cone condition of $F$ holds, namely,
  \begin{equation}\label{equation 2-1}
    \left\|F(x)-F(\tilde{x})-F'(x)(x-\tilde{x})\right\|\leq
    \eta\left\|F(x)-F(\tilde{x})\right\|,~ \forall x, \tilde{x}\in B_{4\rho}(x_0),
  \end{equation}
  where $0<\eta<1$.
  \item[(3)]
  Fr\'{e}chet derivative $F'(\cdot)$ is locally uniformly nonzero and bounded, i.e.,
  \[0<\left\|F'(x)\right\|\leq c_F,~~~ \forall x\in B_{4\rho}(x_0),\]
  where $c_F>0$.
\end{description}
\end{assumption}

\begin{remark}\label{remark 2-1}
Assume that $\eta=0$, we can get the operator F is linear. Since we are concerned with the non-linear case, we will thus occasionally ignore the case where $\eta=0$ without any further remarks.

 The validity of local tangential cone condition indicates that the operator $F'$ fulfills $F'(x)=0$ for some $x\in B_{4\rho}(x_0)$ if and only if $F$ is constant in $B_{4\rho}(x_0)$ (see e.g. \cite[Proposition 1.12]{22}). The case where $F$ is constant in $B_{4\rho}(x_0)$ is not interested, we thus postulate $F'(x)\neq0$ for all $x\in B_{4\rho}(x_0)$.
\end{remark}

\begin{lemma} \label{lemma 2-1}
(\cite[Proposition 2.1]{2}). Let $\rho, \varepsilon>0$ be such that \[\left\|F(x)-F(\tilde{x})-F'(x)(x-\tilde{x})\right\|\leq
    c(x,\tilde{x})\left\|F(x)-F(\tilde{x})\right\|, \forall x, \tilde{x}\in B_{\rho}(x_0)\subset\mathcal{D}(F)\]
for some $c(x,\tilde{x})\geq0,$ and $c(x,\tilde{x})<1$ if $\|x-\tilde{x}\|\leq\varepsilon$. Assume that $F(x)=y$ is solvable in $B_{\rho}(x_0)$, then a unique $x_0$-minimum-norm solution exists which is characterized as the solution $x^\dagger$ of $F(x)=y$ in $B_{\rho}(x_0)$ satisfying
\[x^\dagger-x_0\in\mathcal{N}\left(F'(x^\dagger)\right)^\perp.\]
\end{lemma}

Now, we give some definitions and their preliminary results about metric projection, which are crucial for our further analysis.

 \begin{definition} \label{definition 2-2}
 (\cite[Lemma 2.1]{22}).
 The metric projection $P_C(x)$ satisfies
 \begin{equation}\label{equation 2-2}
 \|x-P_C(x)\|^2=\mathop{\min}\limits_{z\in C}\|x-z\|^2,
 \end{equation}
 where $C\subset \mathcal{X}$ is a nonempty closed convex set which can guarantee $P_C(x)$ is the unique element of $x\in \mathcal{X}$.
 For later convenience, we use the square of the distance. The metric projection onto a convex set $C$ fulfills a descent property which reads
 \begin{equation}\label{equation 2-3}
 \|P_C(x)-z\|^2\leq\|x-z\|^2-\|P_C(x)-x\|^2
 \end{equation}
 for all $z\in C$.
 \end{definition}

 \begin{definition}\label{definition 2-3}
 Define the $\emph{hyperplane}$
 \[H(u,\alpha):=\{x\in\mathcal{X}:\langle u,x \rangle=\alpha\},\]
 where $u\in\mathcal{X}\backslash\{0\}$, $\alpha,\xi\in R$ with $\xi\geq0$ . The halfspace \[H_\leq(u,\alpha):=\{x\in\mathcal{X}:\langle u,x \rangle\leq\alpha\}.\]
 We also can define $H_\geq(u,\alpha)$, $H_<(u,\alpha)$ and $H_>(u,\alpha)$ in the same way. Then, the stripe
 \[H(u,\alpha,\xi):=\{x\in\mathcal{X}:|\langle u,x \rangle-\alpha|\leq\xi\}.\]
 Note that we can get $H(u,\alpha,\xi)=H_\leq(u,\alpha+\xi)\cap H_\geq(u,\alpha-\xi)$ and $H(u,\alpha,0)=H(u,\alpha)$. In addition, $H(u,\alpha)$, $H_\leq(u,\alpha)$, $H_\geq(u,\alpha)$ as well as $H(u,\alpha,\xi)$ are nonempty, closed, convex subsets of $\mathcal{X}$. Meanwhile, $H_<(u,\alpha)$, $H_>(u,\alpha)$ are open subsets of $\mathcal{X}$.
 \end{definition}

 \begin{remark}\label{remark 2-2}
 The metric projection $P_{H(u,\alpha)}(x)$ in the Hilbert space setting corresponds to an orthogonal projection, which reads
 \begin{equation}\label{equation 2-4}
 P_{H(u,\alpha)}(x)=x-\frac{\langle u,x \rangle-\alpha}{\|u\|^2}u.
 \end{equation}
 \end{remark}

  \noindent The following statements can be proved in \cite{15, 16, 23}, which are very helpful to our method.

\begin{prop}\label{proposition 2-4}
The following statements hold.
\begin{description}
  \item[(1)] Let $H(u_i,\alpha_i), \;i=1,2,...,N,$ be hyperplanes with nonempty intersection $H:=\mathop{\bigcap}\limits_{i=1,2,...,N}H(u_i,\alpha_i)$. The projection of x onto H is calculated by
      \[P_H(x)=x-\sum\limits_{i=1}^N \tilde{t}_i u_i,\]
  where $\tilde{t}=(\tilde{t}_1,...,\tilde{t}_N)\in\mathbb{R}^N$ is the solution of the following convex optimization problem
  \begin{equation}\label{equation 2-5}
   \mathop{\min}\limits_{t\in R^N}h(t):=\frac{1}{2}\left\|x-\sum\limits_{i=1}^N t_iu_i\right\|^2+\sum\limits_{i=1}^N t_i \alpha_i.
   \end{equation}
   The corresponding partial derivatives are
  \begin{equation}\label{equation 2-6}
   \frac{\partial}{\partial t_j}h(t)=-\left\langle u_j,x-\sum\limits_{i=1}^N t_iu_i \right\rangle+\alpha_j.
   \end{equation}
  If vectors $u_i,\; i=1,2,...,N$, are linearly independent, we can get that $h$ is strictly convex and $\tilde{t}$ is unique.
  \item[(2)] If $x\in H_>(u,\alpha)$, the projection of x onto $ H_\leq(u,\alpha)$ is obtained as follows
      \[P_{H_\leq(u,\alpha)}(x)=P_{H(u,\alpha)}(x)=x-t_+u,\]
      where\[t_+=\frac{\langle u,x \rangle-\alpha}{\|u\|^2}>0.\]
  \item[(3)] The projection of $x\in\mathcal{X}$ onto a stripe $H(u,\alpha,\xi)$ can be calculated as
      \[P_{H(u,\alpha,\xi)}(x)=
      \left\{
    \begin{array}{ll}
      P_{H_\leq(u,\alpha+\xi)}(x),~~\;\; x\in H_>(u,\alpha+\xi); \\
      \quad \quad x,~~~~~~~~~~\;\quad x\in H(u,\alpha,\xi); \\
      P_{H_\geq(u,\alpha-\xi)}(x),~~\;\; x\in H_<(u,\alpha-\xi).
    \end{array}
  \right.
  \]
\end{description}
\end{prop}

\begin{definition}\label{definition 2-5}
Since the projection point is not easy to find in the proposed algorithm, we specify the order of the projections:
\begin{enumerate}
  \item First, solve for the metric projection point $p_{k_1}:=P_{H_{k,k}}(z_k)$ from $z_k$ to $H_{k,k}$.
  \item And then, get  $p_{k_2}:=P_{H_{k,k_1}\mathop{\bigcap}H_{k,k_2}}(P_{H_{k,k}}(z_k))$ in the same way. Meanwhile, we can get $p_{k_3}:=P_{H_{k,k_1}\mathop{\bigcap}H_{k,k_2}\mathop{\bigcap}H_{k,k_3}}(p_{k_2})$.
  \item We end up with $P_{\mathop{\mathop{\bigcap}\limits_{i\in I_{k}}H_{k,i}}}(p_{k_{s-1}})$ as $x_{k+1}$ just like the previous definitions.
\end{enumerate}

We sort the elements of the finite set $I_k$ that contains $k$, such as $k=k_1>k_2>\cdots>k_s$. And, $H_{k,i}:=H(u_{k,i},\alpha_{k,i},\xi_{k,i})$, where $i\in I_k$. Obviously, our definition of $H_{k,i}$ and their finite intersection are closed convex sets. According to the optimum approximation theorem, we know that there is a unique projection point onto a closed convex subset in Hilbert spaces.
\end{definition}

\section{The TGSS method}\label{section-3}

 \noindent In this section, we describe in detail the TGSS method with exact data case and noisy data case respectively. Then, we formulate the TGSS scheme with two search directions.

\subsection{TGSS with exact data}\label{subsection 3-1}
\noindent Firstly, we put forward TGSS method for nonlinear operators in the case of exact data. We give the following definition of iteration.
\begin{definition}\label{definition 3-1}
At iteration $k\in \mathbb{N}$, choose a finite index set $I_k$  includes $k$. Let \begin{equation}\label{equation 3-1}
z_{i}=x_{i}+\lambda_{i}(x_{i}-x_{i-1}),
\end{equation}
where
$\lambda_{0}=0,~ 0\leq \lambda_{i}\leq 1,~\forall i\in I_k\subset \mathbb{N},$ and define
\[H_{k,i}:=H(u_{k,i},\alpha_{k,i},\xi_{k,i})\]
with
\begin{equation}\label{equation 3-2}
\begin{array}{ll}
 w_{k,i}:=F(z_{i})-y,\\
 u_{k,i}:=F'(z_{i})^{*}w_{k,i},\\
 \alpha_{k,i}:=\langle u_{k,i},z_{i}\rangle-\langle w_{k,i},r_{i}\rangle,\\
 \xi_{k,i}:=\eta\| w_{k,i}\|\| r_{i}\|,
\end{array}
\end{equation}
where $r_{i}=F(z_{i})-y$ is the residual term, which in this case is equal to $w_{k,i}$.
\end{definition}

\begin{alg}\label{algorithm 3-1}
\rm(\textbf{TGSS iteration for exact data with multiple search directions.})

\textbf{Given:}~Exact data $y$, initial choice $x_0=x_{-1}\in\mathcal{X}$ and the parameter $\eta$.

\textbf{For}~$k\in \mathbb{N}$.

~~~~~\textbf{If}~$F(x_{k})\neq y$.

~~~~~~~~~(1) Choose a index set $I_{k}\subset \{k+K+1, \ldots, k\}\cap \mathbb{N}$ with $K\geq1$.

~~~~~~~~~~(2) Compute $z_k$, the search directions $u_{k,i}$,  and the parameters $\alpha_{k,i}$, $\xi_{k,i}$ defined in Definition \ref{definition 3-1}.

~~~~~~~~~(3) Update $x_{k+1}$ by
\begin{equation}\label{equation 3-3}
 x_{k+1}=z_{k}-\sum\limits_{i\in I_{k}}t_{k,i}u_{k,i}.
\end{equation}
~~~~~~~~~~~~~~~
Find $t_k:=(t_{k,i})_{i\in I_k}$ by Definition \ref{definition 2-5} such that
\begin{equation}\label{equation 3-4}
x_{k+1}\in \bigcap_{i\in I_{k}}H_{k,i}.
\end{equation}
~~~~~~~~~~~~~~~
Output~$x_{k+1}$.

~~~~~~~~~
Set $k=k+1$.

~~~~~\textbf{Else if}~$F(x_{k})=y$.

~~~~~~~~~
Output~$x_{k}$.

~~~~~~~~~
\textbf{break}

~~~~~\textbf{End If}

\textbf{End For}
\end{alg}


\begin{prop}\label{proposition 3-2}
For any $k\in \mathbb{N}, i\in I_k$, the solution set $M_{F(x)=y}$ fulfills
\begin{equation}\label{equation 3-5}
  M_{F(x)=y}\subset H_{k,i},
\end{equation}
where $u_{k,i},\alpha_{k,i}$ and $\xi_{k,i}$ are chosen as in (\ref{equation 3-2}).
\begin{proof} Let $z\in M_{F(x)=y}$, we have
\[\langle u_{k,i},z\rangle-\alpha_{k,i}=\langle w_{k,i},F'(z_{i})(z-z_{i})+F(z_{i})-y\rangle.\]
Due to $F(z)=y$ there hold
\begin{equation*}
\begin{array}{ll}
\left|\langle w_{k,i},F'(z_{i})(z-z_{i})+F(z_{i})-F(z)\rangle\right|\leq\|w_{k,i}\|\cdot\|F(z_{i})-F(z)-F'(z_{i})(z_{i}-z)\|\\
~~~~~~~~~~~~~~~~~~~~~~~~~~~~~~~~~~~~~~~~~~~~~~~~~\leq\eta\|w_{k,i}\|\cdot\|F(z_{i})-y\|,
    \end{array}
 \end{equation*}
and using $r_{i}=F(z_{i})-y$ we have $z\in H_{k,i}$.
\end{proof}
\end{prop}

\begin{prop}\label{proposition 3-3}
 Let $\{x_k\}_{k\in \mathbb{N}},\{z_k\}_{k\in N}$ be iterative sequences generated by Algorithm \ref{algorithm 3-1}.  We have
\begin{equation}\label{equation 3-6}
 z_{k}\in H_>(u_{k,k},\alpha_{k,k}+\xi_{k,k}).
\end{equation}
 By projecting $z_k$ first onto $H(u_{k,k},\alpha_{k,k},\xi_{k,k})$ we obtain the descent property
\begin{equation}\label{equation 3-7}
  \left\|z-x_{k+1}\right\|^2\leq\left\|z-z_{k}\right\|^2-\frac{(1-\eta)^2\left\|r_k\right\|^4}{\left\|u_{k,k}\right\|^2},
  \end{equation}
  for $z\in M_{F(x)=y}.$
  \begin{proof}
The first estimate is due to an adequate choice of the stripes. According to our defination of $w_{k,k}=F(z_{k})-y=r_{k}$ we have
\begin{equation*}
\begin{array}{ll}
\alpha_{k,k}:=\langle u_{k,k},z_{k}\rangle-\|r_{k}\|^2,\\
\xi_{k,k}:=\eta\| r_{k}\|^2,
\end{array}
\end{equation*}
and thus $\langle u_{k,k},z_{k}\rangle-\alpha_{k,k}=\|r_{k}\|^2>\xi_{k,k}$ as $0<\eta<1$. According to the Proposition \ref{proposition 2-4} and Definition \ref{definition 2-5}
we have
\begin{equation*}
\begin{array}{ll}
  \left\|z-x_{k+1}\right\|^2\leq\left\|z-P_{H(u_{k,k},\alpha_{k,k},\xi_{k,k})}
  (z_{k})\right\|^2\\
  ~~~~~~~~~~~~~~~~=\left\|z-z_k+\frac{\left\langle u_{k,k},z_k \right\rangle-(\alpha_{k,k}+\xi_{k,k})}{\left\|u_{k,k}\right\|^2}u_{k,k}\right\|^2
  =\left\|z-z_k+\frac{(1-\eta)\left\|F(z_k)-y\right\|^2}{\left\|u_{k,k}\right\|^2}u_{k,k}\right\|^2\\
  ~~~~~~~~~~~~~~~~=\left\|z-z_k\right\|^2+2\frac{(1-\eta)\left\|F(z_k)-y\right\|^2}{\left\|u_{k,k}\right\|^2}
  \left\langle u_{k,k},z-z_k \right\rangle+\frac{(1-\eta)^2\left\|F(z_k)-y\right\|^4}{\left\|u_{k,k}\right\|^2}\\
  ~~~~~~~~~~~~~~~~=\left\|z-z_k\right\|^2+2\frac{(1-\eta)\left\|r_k\right\|^2}{\left\|u_{k,k}\right\|^2}
  \left\langle r_{k},F'(z_{k})(z-z_k)+r_k-r_k \right\rangle+\frac{(1-\eta)^2\left\|r_k\right\|^4}{\left\|u_{k,k}\right\|^2}\\
  ~~~~~~~~~~~~~~~~\leq\left\|z-z_k\right\|^2-\frac{(1-\eta)^2\left\|r_k\right\|^4}{\left\|u_{k,k}\right\|^2}.
  \end{array}
\end{equation*}
The descent property (the second estimate) is  obtained.
  \end{proof}
\end{prop}

\subsection{TGSS with noisy data}\label{subsection 3-2}

\noindent We have discussed the case of exact data before which is an idealized case in practice. Now we will extend our discussion to noisy data case in this subsection.

\begin{definition}\label{definition 3-4}
At iteration $k\in \mathbb{N}$, choose a finite index set $I_k^\delta$ includes $k$. Let
 \begin{equation}\label{equation 3-8}
 z_{i}^\delta=x_{i}^\delta+\lambda_{i}^\delta(x_{i}^\delta-x_{i-1}^\delta),
 \end{equation}
where
$\lambda_{0}^\delta=0,\;0\leq \lambda_{i}^\delta\leq 1,\;\forall i\in I_k^\delta\subset \mathbb{N}$.
Define the stripes
\[H_{k,i}^\delta:=H(u_{k,i}^\delta,\alpha_{k,i}^\delta,\xi_{k,i}^\delta)\]
with
\begin{equation}\label{equation 3-9}
\begin{array}{ll}
 w_{k,i}^\delta:=F(z_{i}^\delta)-y^\delta,\\
 u_{k,i}^\delta:=F'(z_{i}^\delta)^{*}w_{k,i}^\delta,\\
 \alpha_{k,i}^\delta:=\langle u_{k,i}^\delta,z_{i}^\delta\rangle-\langle w_{k,i}^\delta,r_{i}^\delta\rangle,\\
 \xi_{k,i}^\delta:=(\delta+\eta(\| r_{i}^\delta\|+\delta))\| w_{k,i}^\delta\|,
\end{array}
\end{equation}
where $r_{i}^\delta=F(z_{i}^\delta)-y^\delta$ is the residual term, which in this case is equal to $w_{k,i}^\delta$.
\end{definition}

\noindent In this paper, we will use the $Morozov~discrepency~principle$ with respect to $z_{k}^\delta$, i.e., stop the iteration after $k_*$ steps, where
\begin{equation}\label{equation 3-10}
k_*:=k_*(\delta,y^\delta):=\min\{k\in \mathbb{N}:\|r_k^\delta\|\leq\tau\delta\}.
\end{equation}

\begin{alg}\label{algorithm 3-2}
\rm(\textbf{TGSS iteration for noisy data with multiple search directions.})

\textbf{Given:}~Noisy data $y^\delta$, parameter $\delta$ and initial choice $x_0^\delta=x_{-1}^\delta=x_0\in\mathcal{X}$ and the parameter $\eta$.

\textbf{For}~$k\in \mathbb{N}$.

~~~~~\textbf{If}~stopping criterion (\ref{equation 3-10}) is not satisfied.

~~~~~~~~~(1) Choose a index set $I_{k}^\delta\subset \{k+K+1, \ldots, k\}\cap \mathbb{N}$ with $K\geq1$.

~~~~~~~~~~(2) Compute $z_k^\delta$, the search directions $u_{k,i}^\delta$,  and the parameters $\alpha_{k,i}^\delta$, $\xi_{k,i}^\delta$ defined in Definition \ref{definition 3-4}.

~~~~~~~~~(3) Update $x_{k+1}^\delta$ by
\begin{equation}\label{equation 3-11}
 x_{k+1}^\delta=z_{k}^\delta-\sum\limits_{i\in I_{k}^\delta}t_{k,i}^\delta u_{k,i}^\delta.
\end{equation}
~~~~~~~~~~~~~~~
Find $t_k^\delta:=(t_{k,i}^\delta)_{i\in I_k^\delta}$ by Definition \ref{definition 2-5} such that
\begin{equation}\label{equation 3-12}
x_{k+1}^\delta\in \bigcap_{i\in I_{k}^\delta}H_{k,i}^\delta.
\end{equation}
~~~~~~~~~~~~~~~
Output~$x_{k+1}^\delta$.

~~~~~~~~~
Set $k=k+1$.

~~~~~\textbf{Else if}~$k$ satisfies (\ref{equation 3-10}).

~~~~~~~~~
Output~$x_{k}$.

~~~~~~~~~
Set $k=k_*$.

~~~~~~~~~
\textbf{break}

~~~~~\textbf{End If}

\textbf{End For}
\end{alg}

\noindent We obtain the analogous statement as in the noise-free case.

\begin{prop}\label{proposition 3-5}
For any $k\in \mathbb{N}, i\in I_k^\delta$, the solution set $M_{F(x)=y}$ fulfills
\[M_{F(x)=y}\subset H_{k,i}^\delta,\]
 where $u_{k,i}^\delta,\alpha_{k,i}^\delta$ and $\xi_{k,i}^\delta$ are chosen as in (\ref{equation 3-9}).

\begin{proof} Let $z\in M_{F(x)=y}.$ We then have
\begin{equation*}
\begin{array}{ll}
|\langle u_{k,i}^\delta,z\rangle-\alpha_{k,i}^\delta|=|\langle w_{k,i}^\delta,F'(z_{i}^\delta)(z-z_{i}^\delta)+F(z_{i}^\delta)-F(z)+F(z)-y^\delta\rangle|\\
~~~~~~~~~~~~~~~~~~~~\leq\|w_{k,i}^\delta\|\cdot(\|F(z_{i}^\delta)-F(z)-F'(z_{i}^\delta)(z_{i}^\delta-z)\|+\|y^\delta-y\|)\\
~~~~~~~~~~~~~~~~~~~~\leq\|w_{k,i}^\delta\|\cdot(\eta\|F(z_{i}^\delta)-y^\delta+y^\delta-F(z)\|+\delta)\\
~~~~~~~~~~~~~~~~~~~~\leq\|w_{k,i}^\delta\|\cdot(\eta(\| r_{i}^\delta\|+\delta)+\delta),
\end{array}
\end{equation*}
i.e., $z\in H_{k,i}^\delta$.
\end{proof}
\end{prop}

 \begin{prop}\label{proposition 3-6}
 Let $\{x_k^\delta\}_{k\in \mathbb{N}},\{z_k^\delta\}_{k\in N}$ be iterative sequences generated by Algorithm \ref{algorithm 3-2}. As long as $\|r_{k}^\delta\|>\frac{1+\eta}{1-\eta}\delta$, we have
\begin{equation}\label{equation 3-13}
 z_k^\delta\in H_>(u_{k,k}^\delta,\alpha_{k,k}^\delta+\xi_{k,k}^\delta),
\end{equation}
where $u_{k,k}^\delta, \alpha_{k,k}^\delta$ and $\xi_{k,k}^\delta$ are chosen as in (\ref{equation 3-9}). By projecting $z_k^\delta$ first onto $H(u_{k,k}^\delta,\alpha_{k,k}^\delta,\xi_{k,k}^\delta)$ we get similar decreasing property to Proposition \ref{proposition 3-3} (the second estimate) which we discuss in details later.
  \begin{proof}
The proof we refer to Proposition \ref{proposition 3-3} and Definition \ref{definition 3-4}.
  \end{proof}
 \end{prop}

In the following subsection, we want to take a look at an important special case of Algorithm \ref{algorithm 3-2} which acquires a better understanding of the structure of the TGSS method.

\subsection{TGSS with two search directions}\label{subsection 3-3}

\noindent We want to summarize a fast way to compute $x_{k+1}^\delta$ according to Algorithm \ref{algorithm 3-2}, using only two search directions, i.e., $I_{k}^\delta$ has two elements ($N=2$). This method has been suggested and analyzed by Sch\"{o}pfer and Schuster in \cite{23} and has been successfully implemented for the numerical solution of an integral equation of the first kind.

The following algorithm is a special case of Algorithm \ref{algorithm 3-2}, where we have chosen $I_{0}^\delta=\{0\}$ and $I_{k}^\delta=\{k-1,k\}$ for all $k\geq1$. For convenience, we skip the first index $k$ in the subscript of the functions and parameters we are dealing with.

\begin{definition}\label{definition 3-7}
In the first step (k=0), we choose $u_0^\delta$ as the search direction. At iteration $k\in \mathbb{N}$, where $k\geq 1$ and $\mathbb{N}$ is the set of all natural numbers, choose index set $I_k^\delta:=\{k,k-1\}$, i.e., the two search directions $\{u_k^\delta,u_{k-1}^\delta\}$, where
\begin{equation}\label{equation 3-14}
\begin{array}{l}
w_{k}^\delta:=F(z_{k}^\delta)-y^\delta,\\
u_k^\delta:=F'(z_{k}^\delta)^{*}w_{k}^\delta.
\end{array}
\end{equation}
Let $H_{-1}^\delta:=\mathcal{X}$, and for $k\in \mathbb{N}$, define the stripes
\[H_{k}^\delta:=H(u_{k}^\delta,\alpha_{k}^\delta,\xi_{k}^\delta)\]
with
\begin{equation}\label{equation 3-15}
\begin{array}{l}
 \alpha_{k}^\delta:=\langle u_{k}^\delta,z_{k}^\delta\rangle-\langle w_{k}^\delta,r_{k}^\delta\rangle,\\
 \xi_{k}^\delta:=(\delta+\eta(\| r_{k}^\delta\|+\delta))\| w_{k}^\delta\|.
\end{array}
\end{equation}

\noindent Where $r_{k}^\delta=F(z_{k}^\delta)-y^\delta$, is the residual term of iteration which in this case is equal to $w_{k}^\delta$. Choose  the same stop rule as in Algorithm \ref{algorithm 3-2}, where the constant holds
\begin{equation}\label{equation 3-16}
\tau>\frac{1+\eta}{1-\eta}>1.
\end{equation}
\end{definition}

\noindent If $\left\|r_{k}^\delta\right\|>\tau\delta$, it follows Proposition \ref{proposition 3-6} that
\begin{equation}\label{equation 3-17}
z_k^\delta\in H_>(u_{k}^\delta,\alpha_{k}^\delta,\xi_{k}^\delta).
\end{equation}
Then calculate the iterate point $x_{k+1}^\delta$ by the following steps.
\begin{description}
  \item[(i)] Compute
  \[\tilde{x}_{k+1}^\delta:=P_{H(u_{k}^\delta,\alpha_{k}^\delta+\xi_{k}^\delta)}(z_k^\delta)=z_k^\delta-\frac{\langle u_{k}^\delta,z_k^\delta \rangle-(\alpha_{k}^\delta+\xi_{k}^\delta)}{\|u_{k}^\delta\|^2}u_{k}^\delta.\]
  Thus, for all $z\in M_{F(x)=y}$, we have the descent property
  \begin{equation}\label{equation 3-18}
  \left\|z-\tilde{x}_{k+1}^\delta\right\|^2\leq\left\|z-z_k^\delta\right\|^2-
  \left(\frac{\|r_{k}^\delta\|\left(\|r_{k}^\delta\|-\delta-\eta(\| r_{k}^\delta\|+\delta)\right)}{\|u_{k}^\delta\|}\right)^2.
  \end{equation}
  If $\tilde{x}_{k+1}^\delta\in H_{k-1}^\delta$, i.e., $\tilde{x}_{k+1}^\delta=P_{H_{k}^\delta\cap H_{k-1}^\delta}(z_{k}^\delta)$, the calculation is completed. Otherwise, turn to step (ii).
  \item[(ii)] Decide whether $\tilde{x}_{k+1}^\delta\in H_>(u_{k-1}^\delta,\alpha_{k-1}^\delta+\xi_{k-1}^\delta)$ or $\tilde{x}_{k+1}^\delta\in H_<(u_{k-1}^\delta,\alpha_{k-1}^\delta-\xi_{k-1}^\delta)$.
  Then calculate
  \[x_{k+1}^\delta:=P_{H(u_{k}^\delta,\alpha_{k}^\delta+\xi_{k}^\delta)\cap H(u_{k-1}^\delta,\alpha_{k-1}^\delta\pm\xi_{k-1}^\delta)}(\tilde{x}_{k+1}^\delta),\]
  i.e., determine $x_{k+1}^\delta=\tilde{x}_{k+1}^\delta-\tilde{t}_{k}^\delta u_{k}^\delta-\tilde{t}_{k-1}^\delta u_{k-1}^\delta$ such that $\tilde{t}=(\tilde{t}_{k}^\delta ,\tilde{t}_{k-1}^\delta )$ minimizes the optimization function
 \begin{equation}\label{equation 3-19}
  h_2(t_1,t_2):=\frac{1}{2}\left\|\tilde{x}_{k+1}^\delta-t_1 u_{k}^\delta-t_2 u_{k-1}^\delta \right\|^2+t_1(\alpha_{k}^\delta+\xi_{k}^\delta)+t_2(\alpha_{k-1}^\delta\pm\xi_{k-1}^\delta).
  \end{equation}
  Then $x_{k+1}^\delta\in H_{k}^\delta\cap H_{k-1}^\delta$ and for all $z\in M_{F(x)=y}$ we have
  \begin{equation}\label{equation 3-20}
  \left\|z-x_{k+1}^\delta\right\|^2\leq\left\|z-z_k^\delta\right\|^2-S_{k}^\delta,
  \end{equation}
  where
  \[
  S_{k}^\delta:=\left(\frac{\|r_{k}^\delta\|\left(\|r_{k}^\delta\|-\delta-\eta(\| r_{k}^\delta\|+\delta)\right)}{\|u_{k}^\delta\|}\right)^2+\left(\frac{\langle u_{k-1}^\delta,\tilde{x}_{k+1}^\delta \rangle-(\alpha_{k-1}^\delta\pm\xi_{k-1}^\delta)}{\gamma_k\|u_{k-1}^\delta\|}\right)^2
  \]
  and
 \[
  \gamma_k:=\left(1-\left(\frac{\left|\langle u_{k}^\delta,u_{k-1}^\delta \rangle\right|}{\|u_{k}^\delta\|\|u_{k-1}^\delta\|}\right)^2\right)^{\frac{1}{2}}\in \left(0,1\right].
  \]
\end{description}

 \begin{remark}\label{remark 3-1}
 Discussing the TGSS method with two search directions allows a deeper understanding of TGSS method.
 \begin{description}
  \item[(a)]  Although it might not end up with the metric projection point on the intersection, it does guarantee property $\tilde{x}_{k+1}^\delta=P_{H_{k}^\delta\cap H_{k-1}^\delta}(z_{k}^\delta)$ or $x_{k+1}^\delta\in H_{k}^\delta\cap H_{k-1}^\delta$. And according to the uniqueness of the projection path, the uniqueness of the iteration sequence can also be guaranteed.

  \item[(b)] To see that (\ref{equation 3-17}) is valid if $\left\|r_{k}^\delta\right\|>\tau\delta$, we note that (\ref{equation 3-16}) implies
   \[
   \left\|r_{k}^\delta\right\|>\tau\delta>\delta\frac{1+\eta}{1-\eta}.
   \]
   Because of $0\leq\eta<1$ we have
   \[
   \| r_{k}^\delta\|-\eta\| r_{k}^\delta\|-\delta\eta-\delta>0,
   \]
   yielding
   \[
   \alpha_{k}^\delta+\xi_{k}^\delta=\langle u_{k}^\delta,z_k^\delta \rangle-\| r_{k}^\delta\|\cdot(\| r_{k}^\delta\|-\eta\| r_{k}^\delta\|-\delta\eta-\delta)<\langle u_{k}^\delta,z_k^\delta \rangle.
   \]
   Thus $z_k^\delta\in H_>(u_{k}^\delta,\alpha_{k}^\delta+\xi_{k}^\delta)$, i.e., (\ref{equation 3-17}) holds.

  \item[(c)] We try to make the width of each stripe as narrow as possible, which is related to the value of $\eta$ from the tangential cone condition. Similarly, the choice of (\ref{equation 3-16}) for $\tau$  is also depends strongly on the constant $\eta$. The smaller $\eta$ , the better the approximation of F by its linearization, while a large value of $\eta$ , the bigger the corresponding tolerance $\tau\delta$, that is, the residual term (error) after the iteration is stopped becomes larger.
      In addition, the intersection of finding $x_{k+1}^\delta$ becomes correspondingly wider, so that the difference between $x_{k+1}^\delta$ and the true solution will be larger.

  \item[(d)] The discussion of the relationship between  $u_{k}^\delta$ and $u_{k-1}^\delta$ has been stated in \cite{5,15}, and the modification due to step (ii) might be significant, if the search directions $u_{k}^\delta$ and $u_{k-1}^\delta$ fulfill
      \[
     \frac{\left|\langle u_{k}^\delta,u_{k-1}^\delta \rangle\right|}{\|u_{k}^\delta\|\|u_{k-1}^\delta\|}\approx 1,
      \]
      since the case the coefficient $\gamma_k$ is quite small and therefore $S_{k}^\delta$ is large. This can be illustrated by looking at the situation where $u_{k}^\delta\perp u_{k-1}^\delta$: The projection of $z_k^\delta$ onto $H_{k}^\delta$ is already contained in $H_{k-1}^\delta$, such that step (ii) has no effect. This inspires the method of Heber et al. \cite{24}.

  \item[(e)] Algorithm \ref{algorithm 3-2} with two search directions is valuable for an implementation. The search direction $u_{k-1}^\delta$ has already been calculated for the previous iteration and can be reused. Moreover, in each iteration, better $z_{k}^\delta$ is selected to replace $x_{k}^\delta$, and the two-point gradient method is used to optimize the search direction $F'(x_{j}^\delta)^{*}(F(x_{j}^\delta)-y^\delta),~j=\{k,k-1\}$. Hence the costly computations are to choose a suitable $\lambda_{k}^\delta$ and to get the $z_{k}^\delta$, and the determination of $F'(z_{j}^\delta)^{*}r_{k}^\delta$.

  \item[(f)] The relevant conclusions of Algorithm \ref{algorithm 3-2} are also applicable to the noise-free case by setting $\delta=0$, except that the discrepancy principle has to be replaced, such as a maximal number of iterations. The remainder of the proof follows the same lines as in the treatment of noisy case.

\end{description}
\end{remark}


\begin{prop}\label{proposition 3-8}
If $ u_{k}^\delta,u_{k-1}^\delta$ are linearly dependent, step (i) already yields the metric projection of $z_{k}^\delta$ onto $H_{k}^\delta\cap H_{k-1}^\delta$. Furthermore, we have $ H(u_{k}^\delta,\alpha_{k}^\delta+\xi_{k}^\delta)\subset H_{k-1}^\delta$.
\begin{proof}
 Let $u_{k}^\delta=lu_{k-1}^\delta$, $0\neq l\in \mathbb{R}$.
 For $z\in H(u_{k}^\delta,\alpha_{k}^\delta+\xi_{k}^\delta)\subset H_{k-1}^\delta$ we get
\begin{equation*}
\begin{array}{l}
 \alpha_{k-1}^\delta-\xi_{k-1}^\delta\leq\left\langle u_{k-1}^\delta,z \right\rangle\leq\alpha_{k-1}^\delta+\xi_{k-1}^\delta,\\
 l\left\langle u_{k}^\delta,z \right\rangle=l(\alpha_{k}^\delta+\xi_{k}^\delta).
\end{array}
\end{equation*}
 Hence $l(\alpha_{k}^\delta+\xi_{k}^\delta)\in[\alpha_{k-1}^\delta-\xi_{k-1}^\delta,\alpha_{k-1}^\delta+\xi_{k-1}^\delta]$.
 And, for any $x\in H(u_{k}^\delta,\alpha_{k}^\delta+\xi_{k}^\delta)$, we obtain
 \[
 \left\langle u_{k-1}^\delta,x \right\rangle=l\left\langle u_{k}^\delta,x \right\rangle=l(\alpha_{k}^\delta+\xi_{k}^\delta)\in[\alpha_{k-1}^\delta-\xi_{k-1}^\delta,\alpha_{k-1}^\delta+\xi_{k-1}^\delta],
 \]
 showing that $x\in H_{k-1}^\delta$. It follows that $ H(u_{k}^\delta,\alpha_{k}^\delta+\xi_{k}^\delta)\subset H_{k-1}^\delta$ yields the assertion.

 From Definition \ref{definition 2-5} and Algorithm \ref{algorithm 3-2}, we consider the search direction $\sum\limits_{i\in I_{k}^\delta}t_{k,i}^\delta u_{k,i}^\delta$. Inspired by the above discussion and inductive hypothesis method, we figure out that $\{u_{k,i}^\delta\}_{i\in I_{k}^\delta/T_k}$ are linearly independent, where $T_k$ is the index set of $t_{k,i}=0$.  A similar heuristic argument, the details of which we omit. The proposition also applies to Algorithm \ref{algorithm 3-1}.
\end{proof}
\end{prop}

\section{Convergence and regularity analysis}\label{section-4}
\noindent For the analysis of TGSS methods presented in Section \ref{section-3}, using the assumptions which postulated in Section \ref{section-2}, we establish the convergence results of the TGSS method firstly. We then give the regularity analysis of TGSS method with noisy data. Finally, we will discuss DBTS algorithm with the choice of the combination parameter $\lambda_{k}^\delta$ which completes our theory of regularity.

\subsection{Convergence results}\label{subsection 4-1}
\noindent In Section \ref{section-3} we place some restrictions on the combination parameters $\lambda_{k}^\delta$. Minimal requirements on their values are:
\begin{equation}\label{equation 4-1}
\lambda_{0}^\delta=0,~~0\leq\lambda_{k}^\delta\leq1,~~\forall k\in \mathbb{N}.
\end{equation}

\begin{prop}\label{proposition 4-1}
Let $x_{k}^\delta,x_{k-1}^\delta\in B_{\rho}(x_*)$ for $k\in \mathbb{N}$. Assume that
\begin{equation}\label{equation 4-2}
\left\|F(z_{k}^\delta)-y^\delta\right\|>\tau\delta,
\end{equation}
with $\tau$ satisfying
\begin{equation}\label{equation 4-3}
\tau>\frac{1+\eta}{1-\eta}.
\end{equation}
Define
\begin{equation}\label{equation 4-4}
\Delta_{k}:=\left\|x_{k}^\delta-x_*\right\|^2-\left\|x_{k-1}^\delta-x_*\right\|^2,
\end{equation}
 and
 \begin{equation}\label{equation 4-5}
\Psi:=(1-\eta)-\tau^{-1}(1+\eta)>0.
\end{equation}
Then there holds
\begin{equation}\label{equation 4-6}
\Delta_{k+1}^\delta\leq\lambda_{k}^\delta\Delta_{k}^\delta+
\lambda_{k}^\delta(\lambda_{k}^\delta+1)\|x_{k}^\delta-x_{k-1}^\delta\|^2- \frac{\Psi^2\|r_{k}^\delta\|^2}{ c_F^2},
\end{equation}
where $c_F$ is the upper bound of $\|F'(x)\|$.
\begin{proof}
Since $x_{k}^\delta,x_{k-1}^\delta\in B_{\rho}(x_*)$, $x_{*}\in B_{\rho}(x_0)$ and by the triangle inequality, we have $x_{k}^\delta,x_{k-1}^\delta\in B_{2\rho}(x_0)$. Together with $\lambda_{k}^\delta\leq1$, we deduce
\begin{equation*}
\begin{array}{l}
\|z_{k}^\delta-x_0\|\leq\|z_{k}^\delta-x_{k}^\delta\|+\|x_{k}^\delta-x_0\|=\lambda_{k}^\delta\|x_{k}^\delta-x_{k-1}^\delta\|+\|x_{k}^\delta-x_0\|\\
~~~~~~~~~~~~\;\leq\lambda_{k}^\delta\|x_{k}^\delta-x_{*}\|+\lambda_{k}^\delta\|x_{*}-x_{k-1}^\delta\|+\|x_{k}^\delta-x_0\|\\
~~~~~~~~~~~~\;\leq2\lambda_{k}^\delta\rho+2\rho\\
~~~~~~~~~~~~\;\leq4\rho,
\end{array}
\end{equation*}
which indicates that $z_{k}^\delta\in B_{4\rho}(x_0)$. Hence, using Proposition \ref{proposition 3-6} and (\ref{equation 3-18}), we get
\begin{equation*}
\begin{array}{l}
\left\|x_*-x_{k+1}^\delta\right\|^2\leq\left\|x_*- z_{k}^\delta\right\|^2-\left(\frac{\|r_{k}^\delta\|\left(\|r_{k}^\delta\|-\delta-\eta(\| r_{k}^\delta\|+\delta)\right)}{\|u_{k}^\delta\|}\right)^2\\
~~~~~~~~~~~~~~~~~\;\leq\left\|x_*- z_{k}^\delta\right\|^2-\frac{\Psi^2\|r_{k}^\delta\|^2}{ c_F^2}.
\end{array}
\end{equation*}
Now, using the above inequality, we get
\begin{equation*}
\begin{array}{l}
\Delta_{k+1}^\delta=\left\|x_{k+1}^\delta-x_*\right\|^2-\left\|x_{k}^\delta-x_*\right\|^2\\
~~~~~~~\leq\left\|x_*- z_{k}^\delta\right\|^2-\left\|x_{k}^\delta-x_*\right\|^2-\frac{\Psi^2\|r_{k}^\delta\|^2}{ c_F^2}\\
~~~~~~~=2\langle z_{k}^\delta-x_{k}^\delta,x_{k}^\delta-x_* \rangle+\left\|z_{k}^\delta-x_{k}^\delta\right\|^2-\frac{\Psi^2\|r_{k}^\delta\|^2}{ c_F^2}\\
~~~~~~~=-2\lambda_{k}^\delta\langle x_{k-1}^\delta-x_{k}^\delta,x_{k}^\delta-x_* \rangle+(\lambda_{k}^\delta)^2\left\|x_{k}^\delta-x_{k-1}^\delta\right\|^2 -\frac{\Psi^2\|r_{k}^\delta\|^2}{ c_F^2}\\
~~~~~~~=-\lambda_{k}^\delta\left(\left\|x_{k-1}^\delta-x_{k}^\delta+x_{k}^\delta-x_*\right\|^2-\left\|x_{k}^\delta-x_*\right\|^2-\left\|x_{k}^\delta-x_{k-1}^\delta\right\|^2\right)\\
~~~~~~~~~~+(\lambda_{k}^\delta)^2\left\|x_{k}^\delta-x_{k-1}^\delta\right\|^2
-\frac{\Psi^2\|r_{k}^\delta\|^2}{ c_F^2}\\
~~~~~~~=-\lambda_{k}^\delta\left(\left\|x_{k-1}^\delta-x_*\right\|^2-\left\|x_{k}^\delta-x_*\right\|^2\right)+\lambda_{k}^\delta(\lambda_{k}^\delta+1)\left\|x_{k}^\delta-x_{k-1}^\delta\right\|^2-\frac{\Psi^2\|r_{k}^\delta\|^2}{ c_F^2}\\
~~~~~~~=\lambda_{k}^\delta\Delta_{k}^\delta+\lambda_{k}^\delta(\lambda_{k}^\delta+1)\left\|x_{k}^\delta-x_{k-1}^\delta\right\|^2
-\frac{\Psi^2\|r_{k}^\delta\|^2}{ c_F^2},
\end{array}
\end{equation*}
which yields the assertion.
\end{proof}
\end{prop}

\begin{remark}\label{remark 4-1}
 Analogically the fact that $\Delta_{k+1}^\delta\leq0$ for all $k<k_*$, i.e., that $x_{k+1}^\delta$ is a better approximation of $x_*$ than $x_{k}^\delta$ when the discrepancy principle (\ref{equation 3-10}) is not yet satisfied, in the convergence analysis of Landweber iteration. We want our TGSS method to share this property. According to the inequation (\ref{equation 4-6}), we use the $coupling~condition$:
\begin{equation}\label{equation 4-7}
\lambda_{k}^\delta(\lambda_{k}^\delta+1)\|x_{k}^\delta-x_{k-1}^\delta\|^2-\frac{\Psi^2\|r_{k}^\delta\|^2}{ \mu c_F^2}\leq0,
\end{equation}
which holds for all $0\leq k\leq k_*$, and $\mu$ is a constant satisfying $\mu>1$.
Then, we can derive the sufficient condition
\begin{equation}\label{equation 4-8}
\lambda_{k}^\delta(\lambda_{k}^\delta+1)\|x_{k}^\delta-x_{k-1}^\delta\|^2\leq\frac{\left(\Psi\tau\delta\right)^2}{\mu c_F^2},
\end{equation}
 which leads to the choice
\begin{equation}\label{equation 4-9}
 \lambda_{k}^\delta=\min\left\{\sqrt{\frac{(\Psi\tau\delta)^2}{\mu c_F^2\|x_{k}^\delta-x_{k-1}^\delta\|^2}+\frac{1}{4}}-\frac{1}{2},\frac{k}{k+\alpha}\right\},
\end{equation}
where $\alpha\geq3$ is a given number.
Note that in the above formula for $\lambda_{k}^\delta$, inside the $"\min"$ the second argument is taken to be $k/(k+\alpha)$ which is the combination parameter used in Nesterov acceleration scheme. In case the first position is large, this formula may lead to $\lambda_{k}^\delta=k/(k+\alpha)$ and consequently the acceleration effect of Nesterov can be utilized. It also satisfies the requirement (\ref{equation 4-1}). However, when $\delta=0$, we get $\lambda_{k}^\delta=\lambda_{k}^0=0$ by the definition (\ref{equation 4-9}), which leads to the Sequential subspace optimization iteration.

Based on the above facts, we find a sequence $\lambda_{k}^\delta$ by a discrete backtracking search procedure (DBTS), which takes nonzero for $\delta=0$, satisfies the condition (\ref{equation 4-9}), and ensures the acceleration effect at the same time.
\end{remark}

\begin{prop}\label{proposition 4-2}
 Let $k_*=k_*(\delta,y^\delta)$ be chosen according to (\ref{equation 3-10}), and assume that (\ref{equation 4-7}) holds for all $0\leq k<k_*$. Then $x_{k}^\delta\in B_\rho(x_*)$ as in (\ref{equation 3-11}) is well-defined, and we get
\begin{equation}\label{equation 4-10}
\left\|x_{k+1}^\delta-x_*\right\|\leq\left\|x_{k}^\delta-x_*\right\|,~~~~\forall(-1)\leq k<k_*.
\end{equation}
Moreover, $x_{k}^\delta\in B_\rho(x_*)\subset B_{2\rho}(x_0)$ for all $(-1)\leq k<k_*$ and
\begin{equation}\label{equation 4-11}
\sum\limits_{k=0}^{k_*-1}\left\|F(z_{k}^\delta)-y^\delta\right\|^2\leq\frac{ c_F^2}{\overline\mu\Psi^2}\|x_0^\delta-x_*\|^2
\end{equation}
with $\overline\mu:=(\mu-1)/\mu>0.$
\begin{proof}
Taking $k=0$, it follows from (\ref{equation 4-6}) that
\[
\Delta_{1}^\delta\leq\lambda_{0}^\delta\Delta_{0}^\delta+\lambda_{0}^\delta(\lambda_{0}^\delta+1)\left\|x_{0}^\delta-x_{-1}^\delta\right\|^2
-\frac{\Psi^2\|r_{0}^\delta\|^2}{ c_F^2}.
\]
Using (\ref{equation 4-7}) and $\lambda_{0}^\delta=0$, we can deduce
\[
\Delta_{1}^\delta\leq\lambda_{0}^\delta\Delta_{0}^\delta=0,
\]
which implies that $x_{1}^\delta\in B_\rho(x_*)$. We proceed inductively to demonstrate that
\[
\Delta_{k+1}^\delta\leq\lambda_{k}^\delta\Delta_{k}^\delta\leq\cdot\cdot\cdot\leq\lambda_{k}^\delta\lambda_{k-1}^\delta\cdot\cdot\cdot\lambda_{1}^\delta\Delta_{1}^\delta,
\]
where $0\leq\lambda_{k}^\delta\leq1$. As a result, we have
\[\Delta_{k+1}^\delta\leq0~~~~\forall(-1)\leq k<k_*.\]
Then, we obtain (\ref{equation 4-10}), and $x_{k+1}^\delta\in B_\rho(x_*)$, which completes the induction. From $x_*\in B_{\rho}(x_0)$, we derive $x_{k+1}^\delta\in B_{2\rho}(x_0)$ for all $(-1)\leq k<k_*$.

In addition, from (\ref{equation 4-6}), (\ref{equation 4-7}) and (\ref{equation 4-10}) we can deduce that
\[
\Delta_{k+1}^\delta\leq\frac{\Psi^2\|r_{k}^\delta\|^2}{\mu c_F^2}-\frac{\Psi^2\|r_{k}^\delta\|^2}{ c_F^2},
\]
i.e.,
\[
\overline\mu\frac{\Psi^2\|r_{k}^\delta\|^2}{ c_F^2}\leq\left\|x_{k}^\delta-x_*\right\|^2-\left\|x_{k+1}^\delta-x_*\right\|^2,
\]
and hence
\[
\overline\mu\frac{\Psi^2}{ c_F^2}\sum\limits_{k=0}^{k_*-1}\left\|F(z_{k}^\delta)-y^\delta\right\|^2\leq\|x_0^\delta-x_*\|^2-\|x_{k_*}^\delta-x_*\|^2\leq\|x_0^\delta-x_*\|^2.
\]
We have the estimate
\[
\sum\limits_{k=0}^{k_*-1}\left\|F(z_{k}^\delta)-y^\delta\right\|^2\leq\frac{ c_F^2}{\overline\mu\Psi^2}\|x_0^\delta-x_*\|^2,
\]
which yields the assertion.
\end{proof}
\end{prop}

\noindent Under the same assumptions as Proposition \ref{proposition 4-2}, we have the following corollary.

\begin{corollary}\label{corollary 4-3}
 An obvious induction gives the discrepancy principle yields a finite $k_*=k_*(\delta,y^\delta)$ in Algorithm \ref{algorithm 3-2}.
\begin{proof}
Assume that the discrepancy principle is not satisfied for any iteration index $k$, i.e., $\|F(z_{k}^\delta)-y^\delta\|>\tau\delta$ hold for all $k\in \mathbb{N}$. The conclusion of Proposition \ref{proposition 4-2} becomes
\[
\sum\limits_{k=0}^{+\infty}\left\|F(z_{k}^\delta)-y^\delta\right\|^2\leq\frac{ c_F^2}{\overline\mu\Psi^2}\|x_0^\delta-x_*\|^2,
\]
with $\|x_0^\delta-x_*\|\leq\rho$. We thus obtain
\[
\sum\limits_{k=0}^{+\infty}\left\|F(z_{k}^\delta)-y^\delta\right\|^2\leq\frac{ c_F^2}{\overline\mu\Psi^2}\|x_0^\delta-x_*\|^2\leq\frac{ c_F^2}{\overline\mu\Psi^2}\rho^2<+\infty.
\]
Consequently, the sequence $\left\{\|F(z_{k}^\delta)-y^\delta\|\right\}_{k\in \mathbb{N}}$ has to be the null sequence, i.e.,
\begin{equation}\label{equation 4-12}
\lim\limits_{k\rightarrow+\infty}\|F(z_{k}^\delta)-y^\delta\|=0.
\end{equation}
This is a contradiction to our assumption $\tau\delta<\|r_k^\delta\|$ for all $k\in \mathbb{N}$. Then, there must be a finite stopping index $k_*$ fulfilling the discrepancy principle (\ref{equation 3-10}).
\end{proof}
\end{corollary}

\noindent If we are given exact data $y^\delta=y$, i.e., $\delta=0$, then (\ref{equation 4-11}) becomes
\begin{equation}\label{equation 4-13}
\sum\limits_{k=0}^{+\infty}\left\|F(z_{k})-y\right\|^2\leq+\infty,
\end{equation}
as in the case $k_*=+\infty$. Otherwise, if the sum terminates in a finite number of steps, i.e., if $F(z_{k})= y$ for some k, then the iteration is terminated and a solution is found. So, this is not restriction.

\noindent Combining (\ref{equation 4-13}) and the condition (\ref{equation 4-7}), we can obtain that
\begin{equation}\label{equation 4-14}
\sum\limits_{k=0}^{+\infty}\lambda_{k}(\lambda_{k}+1)\|x_{k}-x_{k-1}\|^2<+\infty,
\end{equation}
from which there obviously follows
\begin{equation}\label{equation 4-15}
\lim\limits_{k\rightarrow+\infty}\|F(z_{k})-y\|^2=0,
\end{equation}
and
\begin{equation}\label{equation 4-16}
\lim\limits_{k\rightarrow+\infty}\lambda_{k}(\lambda_{k}+1)\|x_{k}-x_{k-1}\|^2=0.
\end{equation}
If we can prove that $z_k$ converges, then we can get a solution that iteratively converges to $F(x)=y$. To do this, we first show some intermediate results.

\begin{prop}\label{proposition 4-4}
Let $\{x_k\}_{k\in\mathbb{N} }$ be the iterative sequence generated by (\ref{equation 3-3}). We then have $\|x_{k}-x_{*}\|$ converges to a constant, which characterized as $\varepsilon\geq0$. There hold
\[
\lim\limits_{k\rightarrow+\infty}\|z_{k}-x_{*}\|=\varepsilon.
\]
\begin{proof}
From Proposition \ref{proposition 4-2} it follows that the sequence $\{\left\|x_{k}^\delta-x_*\right\|\}_{k\in\mathbb{N}}$ is a bounded monotonically decreasing sequence, then $\left\|x_{k}^\delta-x_*\right\|$ converges to a constant, which characterized as $\varepsilon\geq0$.
According to the definition of (\ref{equation 3-1}), we have the inequality
\[
\|z_{k}-x_{*}\|=\|x_{k}-x_{*}+\lambda_k(x_{k}-x_{k-1})\|\leq\|x_{k}-x_{*}\|+\lambda_k\|x_{k}-x_{k-1}\|,
\]
and from the estimate of (\ref{equation 3-7}), we can obtain
\begin{equation}\label{equation 4-17}
\sqrt{\left\|x_{k+1}-x_{*}\right\|^2+\frac{(1-\eta)^2}{c_F^2}\|r_{k}\|^2}\leq\|z_{k}-x_{*}\|\leq\|x_{k}-x_{*}\|+\lambda_k\|x_{k}-x_{k-1}\|.
\end{equation}
Using (\ref{equation 4-15}) and (\ref{equation 4-16}) we can deduce that
\begin{equation}\label{equation 4-18}
\begin{array}{l}
\lim\limits_{k\rightarrow+\infty}\|F(z_{k})-y\|=0,\\
\lim\limits_{k\rightarrow+\infty}\lambda_{k}\|x_{k}-x_{k-1}\|=0.
\end{array}
\end{equation}
Taking $k\rightarrow+\infty$ in (\ref{equation 4-17}) yields the assertion.
\end{proof}
\end{prop}

\begin{lemma}\label{lemma 4-5}
(\cite[Lemma 2.7]{13}).
Let $x_*\in B_{4\rho}(x_0)$ be a solution of $F(x)=y$ and $x_1,x_2\in B_{4\rho}(x_0)$. Then we get
\begin{equation}\label{equation 4-19}
    \left\|F'(x_1)(x_*-x_2)\right\|\leq
    2(1+\eta)\left\|F(x_1)-y\right\|+(1+\eta)\left\|F(x_2)-y\right\|.
\end{equation}
\end{lemma}

\begin{lemma}\label{lemma 4-6}
 For the iterative sequences generated by Algorithm \ref{algorithm 3-2}, there hold
 \[
 x_k^\delta=x_0+\sum\limits_{i=0}^{k-1}\lambda_{i}^\delta(x_i^\delta-x_{i-1}^\delta)-\sum\limits_{i=0}^{k-1}\sum\limits_{s\in I_i^\delta}t_{i,s}^\delta u_{i,s}^\delta,
 \]
 and
 \[
 x_l^\delta-x_k^\delta=\sum\limits_{i=k}^{l-1}\lambda_{i}^\delta(x_i^\delta-x_{i-1}^\delta)-\sum\limits_{i=k}^{l-1}\sum\limits_{s\in I_i^\delta}t_{i,s}^\delta u_{i,s}^\delta,
 \]
 as well as
 \[
 x_l^\delta-x_{l-1}^\delta=-\sum\limits_{m=0}^{l-2}\sum\limits_{s\in I_m^\delta}(\mathop\prod\limits_{n=m+1}^{l-1}\lambda_{n}^\delta)t_{m,s}^\delta u_{m,s}^\delta-\sum\limits_{s\in I_{l-1}^\delta}t_{l-1,s}^\delta u_{l-1,s}^\delta.
 \]
 \begin{proof}
 Obviously, the first two statements can be obtained from (\ref{equation 1-7}). The third equality can be proved by
induction, which we refer to the proof of Lemma 2.6 in \cite{13}.
\end{proof}
\end{lemma}

\begin{lemma}\label{lemma 4-7}
According to the definitions in Algorithm \ref{algorithm 3-1}, there hold
\begin{equation}\label{equation 4-20}
\lim\limits_{k\rightarrow+\infty}|t_{k,i}|\|F(z_{i})-y\|^2=0.
\end{equation}
\begin{proof}
From Algorithm \ref{algorithm 3-1} and Proposition \ref{proposition 3-8}, we arrive at
\[
\left\|\sum\limits_{i\in I_{k}}t_{k,i}u_{k,i}\right\|=\left\|x_{k+1}-z_k\right\|\leq4\rho,
\]
as well as $\{u_{k,i}\}_{i\in I_k/Tk}$ are linearly independent. Hence, there exist a positive constant $\alpha$ satisfies
\[
\sum\limits_{i\in I_{k}}|t_{k,i}|\left\|u_{k,i}\right\|\leq\alpha\left\|\sum\limits_{i\in I_{k}}t_{k,i}u_{k,i}\right\|.
\]
Therefore, $|t_{k,i}|\left\|u_{k,i}\right\|$ is bounded by a constant $C$ for $i\in I_{k}$, i.e., $|t_{k,i}|\leq\frac{C}{\left\|u_{k,i}\right\|}$.

Next, we turn to the sufficient conditional proof of the conclusion.
It follows from (\ref{equation 3-7}), coupling condition (\ref{equation 4-7}) and Proposition \ref{proposition 4-2}, we get that
\[
\sum\limits_{k=0}^{+\infty}\frac{\left\|F(z_{i})-y\right\|^4}{\|u_{k,i}\|^2}\leq\frac{ 1}{\overline\mu\Psi^2}\|x_0-x_*\|^2.
\]
Thus,
\[
\lim\limits_{k\rightarrow+\infty}\frac{\|F(z_{i})-y\|^2}{\|u_{k,i}\|}=0,
\]
which implies (\ref{equation 4-20}) holds as well.
\end{proof}
\end{lemma}
 Now we state and prove the convergence for TGSS method in the case of exact data, where we apply an additional condition
\begin{equation}\label{equation 4-21}
\sum\limits_{k=0}^{+\infty}\lambda_{k}\|x_{k}-x_{k-1}\|<+\infty.
\end{equation}
Since under the Assumption \ref{assumption 2-1}, $\{\|x_k-x_{k-1}\|\}$ can be bounded by $2\rho$. The sufficient condition of (\ref{equation 4-21}) is given by
\begin{equation}\label{equation 4-22}
\sum\limits_{k=0}^{+\infty}\lambda_{k}<+\infty.
\end{equation}

\begin{thm}\label{theorem 4-8}
Let $N\geq1$ be a fixed integer and $k\in I_k$ for each  iteration $k\in\mathbb{N}$ in Algorithm \ref{algorithm 3-1}.
Assume that $k_*=k_*(0,y)=+\infty$, combination parameters $\lambda_k$ satisfies (\ref{equation 4-1}), (\ref{equation 4-7}) with $\delta=0$ and (\ref{equation 4-21}). Then the iterates $\{z_k\}_{k\in \mathbb{N}}$ generated by Algorithm \ref{algorithm 3-1} converges to a solution $x_*\in B_{\rho}(x_0)\cap\mathcal{D}(F)$ of $F(x)=y$. If in addition $\mathcal{N}\left(F'(x^\dagger)\right)\subset\mathcal{N}\left(F'(x)\right)$ for all $x\in B_{4\rho}(x^\dagger)$, then the sequence $\{z_k\}_{k\in \mathbb{N}}$ converges to $x^\dagger$ as $k\rightarrow+\infty$.
\begin{proof}
Motivated by the proof of Theorem 2.14 from \cite{25}, we will show that $\{z_k\}_{k\in\mathbb{N} }$ is a Cauchy sequence. Define
\begin{equation}\label{equation 4-23}
e_k:=z_k-x_*.
\end{equation}
This is equivalent to show that the sequence $\{e_k\}_{k\in\mathbb{N}}$ is a Cauchy sequence. We have seen the respective proof in Proposition \ref{proposition 4-4} that $\|e_k\|$ is converges to $\varepsilon$.
Given $j\geq k$, we choose some integer $l=l(k,j)\in \{k,k+1,..., j\}$ such that
\begin{equation}\label{equation 4-24}
\|F(z_{l})-y\|\leq\|F(z_{i})-y\|,~~~~\forall k\leq i\leq j.
\end{equation}
There holds
\[
\|e_j-e_k\|\leq\|e_j-e_l\|+\|e_l-e_k\|
\]
as well as
\begin{equation*}
\begin{array}{l}
\|e_j-e_l\|^2=2\langle e_l-e_j,e_l\rangle+\|e_j\|^2-\|e_l\|^2,\\
\|e_l-e_k\|^2=2\langle e_l-e_k,e_l\rangle+\|e_k\|^2-\|e_l\|^2.
\end{array}
\end{equation*}
Let $k\rightarrow+\infty$, the last two terms on the right-hand side converge to $\varepsilon^2-\varepsilon^2=0$. In order to prove $\|e_j-e_k\|\rightarrow0$, we need to show$|\langle e_l-e_j,e_l\rangle|\rightarrow0$ and $|\langle e_l-e_k,e_l\rangle|\rightarrow0$ as $j,k$ tend to infinity.

\noindent We first consider the term
\begin{equation*}
\begin{array}{l}
|\langle e_l-e_k,e_l\rangle|=|\langle z_l-z_k,e_l\rangle|=|\langle x_l-x_k+\lambda_{l}(x_l-x_{l-1})-\lambda_{k}(x_k-x_{k-1}),e_l\rangle|\\
~~~~~~~~~~~~~~~~\;\leq|\langle x_l-x_k,e_l\rangle|+\lambda_{l}|\langle x_l-x_{l-1},e_l\rangle|+\lambda_{k}|\langle x_k-x_{k-1},e_l\rangle|\\
~~~~~~~~~~~~~~~~\;\leq|\langle x_l-x_k,e_l\rangle|+\lambda_{l}\|x_l-x_{l-1}\|\|e_l\|+\lambda_{k}\|x_k-x_{k-1}\|\|e_l\|.
\end{array}
\end{equation*}
Combining with (\ref{equation 4-18}), we get that
\[
\lim\limits_{k\rightarrow+\infty}\left(\lambda_{l}\|x_l-x_{l-1}\|\|e_l\|+\lambda_{k}\|x_k-x_{k-1}\|\|e_l\|\right)=0.
\]
Now consider
\begin{equation*}
\begin{array}{l}
|\langle x_l-x_k,e_l\rangle|=\left|\left\langle \sum\limits_{i=k}^{l-1}\lambda_{i}(x_i-x_{i-1})-\sum\limits_{i=k}^{l-1}\sum\limits_{s\in I_i}t_{i,s}u_{i,s},e_l\right\rangle\right|\\
~~~~~~~~~~~~~~~~~\leq\sum\limits_{i=k}^{l-1}\lambda_{i}\left|\left\langle x_i-x_{i-1},e_l\right\rangle\right|+\sum\limits_{i=k}^{l-1}\sum\limits_{s\in I_i}|t_{i,s}|\left|\left\langle u_{i,s},e_l\right\rangle\right|.
\end{array}
\end{equation*}
Using the boundness of $\|e_l\|$ and (\ref{equation 4-21}), we get the estimation
\[
\sum\limits_{i=k}^{l-1}\lambda_{i}\left|\left\langle x_i-x_{i-1},e_l\right\rangle\right|\leq\sum\limits_{i=k}^{l-1}\lambda_{i}\| x_i-x_{i-1}\|\|e_l\|\leq\sum\limits_{i=k}^{+\infty}\lambda_{i}\| x_i-x_{i-1}\|\|e_l\|<+\infty.
\]
Then it can be deduced that
\begin{equation}\label{equation 4-25}
\lim\limits_{k\rightarrow+\infty}\left(\sum\limits_{i=k}^{l-1}\lambda_{i}\left|\left\langle x_i-x_{i-1},e_l\right\rangle\right|\right)=0.
\end{equation}
According to Lemma \ref{lemma 4-5}, Lemma \ref{lemma 4-6} and Lemma \ref{lemma 4-7}, we have
\begin{equation*}
\begin{array}{l}
\sum\limits_{i=k}^{l-1}\sum\limits_{s\in I_i}|t_{i,s}|\left|\left\langle u_{i,s},e_l\right\rangle\right|=\sum\limits_{i=k}^{l-1}\sum\limits_{s\in I_i}|t_{i,s}|\left|\left\langle F(z_{s})-y,F'(z_s)(z_l-x_*)\right\rangle\right|\\
~~~~~~~~~~~~~~~~~~~~~~~~~~~~\leq\sum\limits_{i=k}^{l-1}\sum\limits_{s\in I_i}|t_{i,s}|\| F(z_{s})-y\|\|F'(z_s)(z_l-x_*)\|\\
~~~~~~~~~~~~~~~~~~~~~~~~~~~~\leq2(1+\eta)\sum\limits_{i=k}^{l-1}\sum\limits_{s\in I_i}|t_{i,s}|\| F(z_{s})-y\|^2\\
~~~~~~~~~~~~~~~~~~~~~~~~~~~~~~~~+(1+\eta)\sum\limits_{i=k}^{l-1}\sum\limits_{s\in I_i}|t_{i,s}|\| F(z_{s})-y\|\| F(z_{l})-y\|\\
~~~~~~~~~~~~~~~~~~~~~~~~~~~~\leq3(1+\eta)\sum\limits_{i=k}^{l-1}\sum\limits_{s\in I_i}|t_{i,s}|\| F(z_{s})-y\|^2.
\end{array}
\end{equation*}
Note that $I_i$ is a finite set. Using (\ref{equation 4-20}) we get that
\begin{equation}\label{equation 4-26}
\lim\limits_{k\rightarrow+\infty}\left(\sum\limits_{i=k}^{l-1}\sum\limits_{s\in I_i}|t_{i,s}|\left|\left\langle u_{i,s},e_l\right\rangle\right|\right)=0.
\end{equation}
According to (\ref{equation 4-25}) and (\ref{equation 4-26}), we arrive at $|\langle x_l-x_k,e_l\rangle|\rightarrow0$, then it follows that $|\langle e_l-e_k,e_l\rangle|\rightarrow0$ as $k$ tends to infinity. Similarly, it can be shown that $|\langle e_l-e_j,e_l\rangle|\rightarrow0$ as $k\rightarrow+\infty$.  Thus, we obtain
\[
\lim\limits_{k\rightarrow+\infty}\|e_j-e_k\|=0,
\]
which deduce that $\{e_k\}_{k\in\mathbb{N}}$ is a Cauchy sequence and the same holds for $\{z_k\}_{k\in\mathbb{N}}$. Then, $\{z_k\}$ converges to a solution $x_*$ of $F(x)=y$ as $\|F(z_{k})-y\|\rightarrow0$ for $k\rightarrow+\infty$.

 Next, we prove the second part of the theorem with the additional condition $\mathcal{N}\left(F'(x^\dagger)\right)\subset\mathcal{N}\left(F'(x)\right)$ for all $x\in B_{4\rho}(x^\dagger)$.  According to the iterates in Algorithm \ref{algorithm 3-1}, we can get
\begin{equation*}
\begin{array}{l}
z_{k+1}-z_k=x_{k+1}+\lambda_{k+1}(x_{k+1}-x_k)-z_k\\
~~~~~~~~~~\;\;\;=-\sum\limits_{i\in I_k}t_{k,i}u_{k,i}+\lambda_{k+1}(x_{k+1}-x_k)\\
~~~~~~~~~~\;\;\;=-(1+\lambda_{k+1})\sum\limits_{i\in I_k}t_{k,i}u_{k,i}+\lambda_{k+1}(z_{k}-x_k)\\
~~~~~~~~~~\;\;\;=-(1+\lambda_{k+1})\sum\limits_{i\in I_k}t_{k,i}u_{k,i}+\lambda_{k+1}\lambda_{k}(x_{k}-x_{k-1}),
\end{array}
\end{equation*}
thus
\[
z_k-z_0=\sum\limits_{s=0}^{k-1}(z_{s+1}-z_{s})=\sum\limits_{s=0}^{k-1}\left(-(1+\lambda_{s+1})\sum\limits_{i\in I_s}t_{s,i}u_{s,i}+\lambda_{s+1}\lambda_{s}(x_{s}-x_{s-1})\right).
\]
Note that $(1+\lambda_{s+1})\sum\limits_{i\in I_s}t_{s,i}u_{s,i}\in \mathcal{R}\left(F'(z_i)^*\right)$ and since
\[
\mathcal{R}\left(F'(z_i)^*\right)\subset\mathcal{N}\left(F'(z_i)\right)^\perp\subset
\mathcal{N}\left(F'(x^\dagger)\right)^\perp~~for~all~ i\in \mathbb{N},
\]
we arrive at
\[
\sum\limits_{s=0}^{k-1}\left(-(1+\lambda_{s+1})\sum\limits_{i\in I_s}t_{s,i}u_{s,i}\right)\in \mathcal{N}\left(F'(x^\dagger)\right)^\perp.
\]
With the conclusion of Lemma \ref{lemma 4-6}, we have
\[
\sum\limits_{s=0}^{k-1}\lambda_{s+1}\lambda_{s}(x_{s}-x_{s-1})\in \mathcal{N}\left(F'(x^\dagger)\right)^\perp,
\]
we conclude that
\[
z_k-z_0\in \mathcal{N}\left(F'(x^\dagger)\right)^\perp~~for~all~ k\in \mathbb{N},
\]
which is also valid for the limit of $z_k$, i.e., $x_*-x_0\in \mathcal{N}\left(F'(x^\dagger)\right)^\perp$.
It follows from Lemma \ref{lemma 2-1} that $x^\dagger$ is the unique solution satisfying the above condition, then we obtain $z_k \rightarrow x^\dagger$ as $k\rightarrow+\infty$.
\end{proof}
\end{thm}

\begin{corollary}\label{corollary 4-9}
Under the assumptions of Theorem \ref{theorem 4-8}, we arrive at $x_k$ converges to $x_*$, where $x_*$ is the limit of $z_k$ as $k\rightarrow+\infty$.
\begin{proof}
According to the iterates in Algorithm \ref{algorithm 3-1} and (\ref{equation 3-7}), we get that
\[
\left\|x_*-x_{k+1}\right\|^2\leq\left\|x_*-z_{k}\right\|^2,
\]
thus
\[
\left\|x_*-x_{k+1}\right\|\rightarrow0~~as ~~k\rightarrow+\infty,
\]
which confirms the statement.
\end{proof}
\end{corollary}

\subsection{Regularity analysis}\label{subsection 4-2}
\noindent We now discuss the Algorithm \ref{algorithm 3-2} be a convergent regularization method with the discrepancy principle (\ref{equation 3-10}) as the stopping rule. To do this, we need to show that $x_k^\delta$ depends continuously on the data $y^\delta$ firstly.


\begin{lemma}\label{lemma 4-10}
 Let $\{x_k\}$ and $\{z_k\}$ be the sequences generated by the Algorithm \ref{algorithm 3-1} with the initial value $x_0^\delta=x_{-1}^\delta=x_0$. Accordingly, $\{x_k^\delta\}$ and $\{z_k^\delta\}$
 generated by the Algorithm \ref{algorithm 3-2} with the noisy data $y^\delta$, where $y^\delta$ satisfies $\|y^\delta-y\|\leq\delta$. Assume that
\begin{equation}\label{equation 4-27}
\lambda_{k}^\delta\rightarrow\lambda_{k}~~~~as~\delta\rightarrow0.
\end{equation}
Then, there hold
\[
x_k^\delta\rightarrow x_k~~and~z_k^\delta\rightarrow z_k~~as~\delta\rightarrow0,
\]
 for a fixed integer $k\in\mathbb{N}$.
\begin{proof}
We prove the assertion by induction, which is closely follows the corresponding proof of \cite[Lemma 4.10]{16}. From $\lambda_0^\delta=0$, we obtain $z_0^\delta=x_{0}^\delta=x_0$. By the definition of (\ref{equation 3-2}), (\ref{equation 3-9}) and the  continuity of $F'(\cdot)$ on $\mathcal{D}(F)$, we arrive at $z_0^\delta,w_{0,i}^\delta,u_{0,i}^\delta,\alpha_{0,i}^\delta$ and $\xi_{0,i}^\delta$ depend continuously on $y^\delta$, where $i\in I_0^\delta=\{0\}$. Since the norm in Hilbert space is continuous,
thus
\[
x_1^\delta=\tilde{x}_1^\delta=z_0^\delta-\frac{\langle u_{0,i}^\delta,z_0^\delta \rangle-(\alpha_{0,i}^\delta+\xi_{0,i}^\delta)}{\|u_{0,i}^\delta\|^2}u_{0,i}^\delta
\]
depends continuously on $y^\delta$, i.e.,
\[x_1^\delta\rightarrow x_1,~~ z_1^\delta\rightarrow z_1,~~as~\delta\rightarrow0,\]
with $\lambda_0^\delta=\lambda_0=0$.

Now, assume that $x_j^\delta$ and $z_j^\delta$ depend
continuously on the data $y^\delta$ for all $j\leq k$.
Since $F$ is continuously  Fr\'{e}chet differentiable on $\mathcal{D}(F)$.
It follows from Proposition \ref{proposition 3-8} that the search directions with the finite set $I_j^\delta/T_k$ are linearly independent. Combining with the strictly convex property of (\ref{equation 2-5}), we have $\{t_{j,i}^\delta\}$ is uniquely determined by (\ref{equation 2-6}), whose coefficients depend
continuously on $y^\delta$. Similar to the discussion for
the case of $k = 1$, we arrive at
\[u_{k,i}^\delta\rightarrow u_{k,i},~~t_{k,i}^\delta\rightarrow t_{k,i},~~ as~\delta\rightarrow0,\]
thus yields
 \[
 x_{k+1}^\delta=z_{k}^\delta-\sum\limits_{i\in I_{k}^\delta}t_{k,i}^\delta u_{k,i}^\delta,
 \]
 depend continuously on $y^\delta$, i.e.,
 \[
 x_{k+1}^\delta\rightarrow x_{k+1},~~ as~\delta\rightarrow0.
 \]
 According to above conclusion and definition (\ref{equation 4-27}), we have
 \[
 z_{k+1}^\delta\rightarrow z_{k+1},~~ as~\delta\rightarrow0,
 \]
 which yields the assertion.
\end{proof}
\end{lemma}

\begin{thm}\label{theorem 4-11}
Let $k_*$ be chosen by the discrepancy principle (\ref{equation 3-10}). Assume that the coupling condition (\ref{equation 4-7}) holds for all $0\leq k\leq k_*$, combination parameters $\lambda_k^\delta$ satisfies (\ref{equation 4-1}), (\ref{equation 4-21}) and (\ref{equation 4-27}). Then $z_{k_*}^\delta$ generated by the Algorithm  \ref{algorithm 3-2} converges to a solution $x_*\in B_\rho(x_0)\cap\mathcal{D}(F)$ of $F(x)=y$ as $\delta\rightarrow0$. In addition, if $\mathcal{N}\left(F'(x^\dagger)\right)\subset\mathcal{N}\left(F'(x)\right)$ for all $x\in B_{4\rho}(x^\dagger)$, then the sequence $z_{k_*}^\delta$ converges to $x^\dagger$ as $\delta\rightarrow0$.
\begin{proof}
 Inspired by the proof of \cite[Theorem 2.10]{13}, let $x_*$ be the limit point of the sequence $\{z_k\}$ given by Algorithm  \ref{algorithm 3-1}. From Corollary \ref{corollary 4-9}, we have $x_k$ also converges to $x_*$. Let $\delta_n\rightarrow0$ as $n\rightarrow+\infty$ as well as $y_n:= y^{\delta_n}$  be a sequence of noisy data satisfies $\|y-y_n\|\leq \delta_n$ and define $k_n:=k_*(\delta_n,y_n)$ be the stopping index determined by discrepancy principle. There have two cases. First, if $k$ is a finite accumulation point of $k_n$, thus we can
assume that $k_n=k$ for all $n\in N$. Then it follows by using the discrepancy principle that
\[
\|F(z_k^{\delta_n})-y_n\|\leq\tau\delta_n.
\]
From Lemma \ref{lemma 4-10}, we get that $z_k^{\delta}$ depends continuously on $y^\delta$  for a fixed $k\in\mathbb{N}$. Then take the limit of $n\rightarrow+\infty$ in the above inequality, which conclude that
\[
z_k^{\delta_n}\rightarrow z_k,~~F(z_k^{\delta_n})\rightarrow F(z_k)=y,~~as~~n\rightarrow+\infty.
\]
This proves that the iteration terminates with $z_k=x_*$ and $z_k^{\delta_n}\rightarrow x_*$ as $\delta_n\rightarrow0$.\\
For the second case, if $k_n\rightarrow+\infty$ as $n\rightarrow+\infty$. For some $k$ with $k_n>k+1$, Proposition \ref{proposition 4-2} and (\ref{equation 4-1}) yield that
\begin{equation*}
\begin{array}{l}
\|z_{k_n}^{\delta_n}-x_*\|=\|x_{k_n}^{\delta_n}+\lambda_{k_n}^{\delta_n}(x_{k_n}^{\delta_n}-x_{k_n-1}^{\delta_n})-x_*\|\\
~~~~~~~~~~~~~~\leq\|x_{k_n}^{\delta_n}-x_*\|+\lambda_{k_n}^{\delta_n}\|x_{k_n}^{\delta_n}-x_*\|+\lambda_{k_n}^{\delta_n}\|x_{k_n-1}^{\delta_n}-x_*\|\\
~~~~~~~~~~~~~~\leq3\|x_{k}^{\delta_n}-x_*\|\\
~~~~~~~~~~~~~~\leq3\|x_{k}^{\delta_n}-x_{k}\|+3\|x_{k}-x_*\|.
\end{array}
\end{equation*}
If we fixed some $\varepsilon>0$, according to Proposition \ref{proposition 4-2} and Proposition \ref{proposition 4-4} that there exist some corresponding fixed $k(\varepsilon)$ such that $\|x_{k}-x_*\|\leq\varepsilon/6$ for all $n>n(\varepsilon,k)$. Furthermore, the iterations in Algorithm \ref{algorithm 3-2} depend continuously on the fixed $k$, thus we can find an $n=n(\varepsilon,k)$ such that $\|x_{k}^{\delta_n}-x_*\|\leq\varepsilon/6$. Hence, it follows that $n$ is chosen sufficiently large enough which satisfies $k_n > k + 1$, we conclude that
\[
\|z_{k_n}^{\delta_n}-x_*\|\leq3\|x_{k}^{\delta_n}-x_{k}\|+3\|x_{k}-x_*\|
\leq3\frac{\varepsilon}{6}+3\frac{\varepsilon}{6}\leq\varepsilon,
\]
and therefore $z_{k_n}^{\delta_n}\rightarrow x_*$ as $n\rightarrow+\infty$, which confirms the first part of the statement.
In addition, if $\mathcal{N}\left(F'(x^\dagger)\right)\subset\mathcal{N}\left(F'(x)\right)$ for all $x\in B_{4\rho}(x^\dagger)$, in this case, $x_*$ can be chosen as $x_*=x^\dagger$. It follows from Theorem \ref{theorem 4-8} that $z_k\rightarrow x^\dagger$, also $x_k\rightarrow x^\dagger$. We can prove the situation by analogy with the previous discussion, which yields the assertion.
\end{proof}
\end{thm}


\subsection{DBTS: the choice of $\lambda_k^\delta$}\label{subsection 4-3}
\noindent In this subsection, we will discuss the choice of the $\lambda_k^\delta$, which leads to a convergent regularization method and also promotes the acceleration effect.

We have already briefly discussed the choice of the combination parameter in Remark \ref{remark 4-1}. However, $\lambda_k^\delta$ decrease to $0$ as $\delta\rightarrow0$ by the choices,  and there is no corresponding two-point gradient acceleration effect. Therefore, in order to satisfy (\ref{equation 4-7}) and (\ref{equation 4-21}) with $\delta\rightarrow0$, we will introduce the discrete backtracking search (DBTS) algorithm proposed in \cite{13}. Define a function $q:\mathbb{R}_0^+\rightarrow\mathbb{R}_0^+$ which satisfies
\begin{equation}\label{equation 4-28}
  q(m_1)\leq q(m_2),~~\forall m_1>m_2,~~\sum\limits_{k=0}^{+\infty}q(k)<+\infty.
\end{equation}
Now, we introduce DBTS algorithm in detail.

\begin{alg}\label{algorithm 4-1}
\rm(\textbf{DBTS algorithm for calculation combination parameters $\lambda_k^\delta$, $k>1$.})

\textbf{Given:}~$x_k^\delta$,$x_{k-1}^\delta$,$\tau,\delta,\Psi,\mu,c_F,\alpha,y^\delta,F,q:
\mathbb{R}_0^+\rightarrow\mathbb{R}_0^+,i_{k-1}\in\mathbb{N},j_{\max}\in\mathbb{N}.$

\textbf{Calculate}~$\|x_k^{\delta}-x_{k-1}^{\delta}\|$ and define
\[{\beta _k}\left( i \right) = \min \left\{ {\frac{{q\left( i \right)}}{{\left\| {x_k^\delta  - x_{k - 1}^\delta } \right\|}},\frac{k}{k+\alpha}} \right\},~~\alpha\geq3.\]

\textbf{For}~{$ j=1,\cdots, j_{\max}$}.

~~~~~~Set $\lambda _k^\delta  = {\beta_k}\left( {{i_{k - 1}} + j} \right)$.

~~~~~~Calculate $z_k^\delta = x_k^\delta+\lambda_k^\delta(x_k^{\delta}-x_{k-1}^{\delta})$.

~~~~~~{\bf If} {$\|y^{\delta}-F(z_k^\delta)\|\leq \tau \delta$}


~~~~~~~~~~${i_k} = {i_{k - 1}} + j,$

~~~~~~~~~~${\bf break}$.

~~~~~~{\bf Else if}~~~{$\lambda _k^\delta \left( {\lambda _k^\delta  + 1} \right){\left\| {x_k^\delta  - x_{k - 1}^\delta } \right\|^2} \le \frac{\Psi^2}{ \mu c_F^2} {\left\| { F\left( {z_k^\delta } \right)-y^\delta} \right\|^2}$}.

~~~~~~~~~~${i_k} = {i_{k - 1}} + j,$

~~~~~~~~~~${\bf break}$,

~~~~~~{\bf Else if}~~~${i_k} = {i_{k - 1}} + {j_{\max }}$.

~~~~~~~~~~Calculate $\lambda_k^\delta$ by (\ref{equation 4-9}).

~~~~~${\bf End~~If}$

${\bf End ~~For}$\\
  { \bf Output:} $\lambda_k^\delta$, $i_k$.
\end{alg}

\begin{remark}\label{remark 4-3}
Comparing with the reference \cite{13}, we have three modifications: The first one is the definition of $\beta _k$  in which we place $\beta _k(i)=k/(k+\alpha)$ instead of $\beta _k(i)=1$, this modification will speed up convergence by making use of the Nesterov acceleration scheme. The second one is in the second "$\textbf{Else if}$" part, where we set $\lambda _k^\delta=0$ by  (\ref{equation 4-9}), the modification can affect convergence speed. The third one is the same position as above, we change the original assignment statement to a judgment statement for a better loop.
\end{remark}

First, we need to verify that the combination parameter $\lambda _k^\delta$ selected from Algorithm \ref{algorithm 4-1} satisfies condition (\ref{equation 4-7}), and (\ref{equation 4-21}). From Algorithm \ref{algorithm 4-1} that when $\|y^{\delta}-F(z_k^\delta)\|> \tau \delta$ happens, i.e., two cases of "$\textbf{Else if}$", which is obvious that $\lambda _k^\delta$ satisfies condition  (\ref{equation 4-7}).

For exact data case, Algorithm \ref{algorithm 4-1} either makes $\lambda _k=\beta_k(i_k)$ or $\lambda _k=0$. It follows from the definition that $i_k\geq i_{k-1}+1$  and $i_k\geq k$. Thus we arrive at
\[
\sum\limits_{k=0}^{+\infty}\lambda_{k}\|x_{k}-x_{k-1}\|\leq\sum\limits_{k=0}^{+\infty}\beta_k(i_k)\|x_{k}-x_{k-1}\|
\leq\sum\limits_{k=0}^{+\infty}q(i_k)\leq\sum\limits_{k=0}^{+\infty}q(k)<+\infty,
\]
therefore condition (\ref{equation 4-21}) holds.

According to the definition of $\lambda_k^\delta$ in Algorithm \ref{algorithm 4-1}, we can not use Theorem \ref{theorem 4-11} to verify the regularization property of the TGSS method,  since the combination parameter is not continuously dependent on $y^\delta$. In fact, $\lambda_k^\delta$ could have multiple values as $\delta\rightarrow0$ which applied to Algorithm \ref{algorithm 3-1} may generate a variety of iterative sequences for noise-free case. For the generated sequences, we will use $\Gamma_{\eta, q}(x_0)$ to denote the set including all the iterative sequences $\{x_k\},\{z_k\}\subset\mathcal{X}$. Without losing generality, the combination parameters $\lambda_k$ in noise-free case are chosen to satisfy
\begin{equation}\label{equation 4-29}
\lambda_{k}(\lambda_{k}+1)\left\|x_{k}-x_{k-1}\right\|^2
\leq\frac{\Psi^2}{ \mu c_F^2}\|F(z_{k})-y\|^2,
\end{equation}
as well as
\begin{equation}\label{equation 4-30}
0\leq\lambda_{k}\leq\min \left\{ {\frac{{q\left( i_k \right)}}{{\left\| {x_k  - x_{k - 1} } \right\|}},\frac{k}{k+\alpha}} \right\},
\end{equation}
where the sequence $\{i_k\}$  of integers in DBTS method satisfying $i_0=0$ and $1\leq i_k-i_{k-1}\leq j_{\max}$ for all $k$.\\
For a given sequence $\{(x_k)\},\{(z_k)\}\in \Gamma_{\eta, q}(x_0)$, we obtain that the corresponding combination parameters $\lambda_{k}$ satisfy condition (\ref{equation 4-7}) and (\ref{equation 4-21}). Thus, it verifies the convergence of $\{x_k\},\{z_k\}$ using Theorem \ref{theorem 4-8}.
The following discussion shows that the corresponding uniform convergence results.

\begin{prop}\label{proposition 4-12}
Assume that all the conditions in Theorem \ref{theorem 4-8} hold. If $\mathcal{N}\left(F'(x^\dagger)\right)\subset\mathcal{N}\left(F'(x)\right)$ for all $x\in B_{4\rho}(x^\dagger)$, then $\forall~\varepsilon>0$, there exist an integer $k(\varepsilon)$, for any sequence $\{(x_k)\}\in \Gamma_{\eta, q}(x_0)$ we get $\|x^\dagger-x_k\|^2<\varepsilon$ for all $k\geq k(\varepsilon)$.
\begin{proof}
This proof closely follows the corresponding proof for \cite[Proposition 3.9]{14}. By contradiction, suppose that the opposite result is true, i.e.,
there exist an $\varepsilon_0>0$  such that
 \begin{equation}\label{equation 4-31}
 \|x^\dagger-x_{k_l}^{(l)}\|^2\geq\varepsilon_0,
 \end{equation}
 for any $l\geq1$ with $\{(x_k^{(l)})\}\in \Gamma_{\eta, q}(x_0)$ and $k_l>l$.
 We will establish the following system. For each $k=0,1,\cdots$, let $\{l_{k,n}\}$ is a strictly increasing subsequence of positive integers and $\{(\hat {x}_k)\}\in\mathcal{X}$. They satisfy the following conditions:
 \begin{description}
   \item[(i)] $\{(\hat {x}_k)\}\in\Gamma_{\eta, q}(x_0)$.
   \item[(ii)] For each fixed $k$, there hold $x_k^{(l_{k,n})}\rightarrow\hat {x}_k$ and $F(x_k^{(l_{k,n})})\rightarrow F(\hat {x}_k)$ as $n\rightarrow+\infty$.
 \end{description}
 Assume that the above system is available, we will deduce a contradiction. According to (i), it follows from Theorem \ref{theorem 4-8} that $\|x^\dagger-\hat {x}_k\|^2\rightarrow0$ as $k\rightarrow+\infty$. Thus we can find a large integer $\hat {k}$ such that
 \[
 \|x^\dagger-\hat {x}_{\hat {k}}\|^2<\varepsilon_0/2.
 \]
 Then we have
\begin{equation}\label{equation 4-32}
\begin{array}{l}
\varepsilon_0/2>\left(\left\|x^\dagger-\hat {x}_{\hat {k}}\right\|^2-\left\|x^\dagger- x_{\hat {k}}^{(l_{\hat {k},n})}\right\|^2\right)+\left\|x^\dagger- x_{\hat {k}}^{(l_{\hat {k},n})}\right\|^2\\
~~~~~~=\left\langle x^\dagger-\hat {x}_{\hat {k}},x_{\hat {k}}^{(l_{\hat {k},n})}-\hat {x}_{\hat {k}}\right\rangle+\left\langle x^\dagger-x_{\hat {k}}^{(l_{\hat {k},n)}},x_{\hat {k}}^{(l_{\hat {k},n})}-\hat {x}_{\hat {k}}\right\rangle+\left\|x^\dagger- x_{\hat {k}}^{(l_{\hat {k},n})}\right\|^2.
\end{array}
\end{equation}
From property (ii), we arrive at
\[
\left\langle x^\dagger-\hat {x}_{\hat {k}},x_{\hat {k}}^{(l_{\hat {k},n})}-\hat {x}_{\hat {k}}\right\rangle+\left\langle x^\dagger-x_{\hat {k}}^{(l_{\hat {k},n)}},x_{\hat {k}}^{(l_{\hat {k},n})}-\hat {x}_{\hat {k}}\right\rangle\rightarrow0~~as~n\rightarrow+\infty.
\]
Then pick $\hat{n}$ with $\hat{l}:=l_{\hat{k},\hat{n}}$ satisfies
\[
\left\langle x^\dagger-\hat {x}_{\hat {k}},x_{\hat {k}}^{(\hat{l})}-\hat {x}_{\hat {k}}\right\rangle+\left\langle x^\dagger-x_{\hat {k}}^{(\hat{l})},x_{\hat {k}}^{(\hat{l})}-\hat {x}_{\hat {k}}\right\rangle\geq-\varepsilon_0/2.
\]
Thus, by (\ref{equation 4-32}) we can obtain
\[
\left\|x^\dagger- x_{\hat {k}}^{(\hat{l})}\right\|^2<\varepsilon_0.
\]
Note that $k_{\hat{l}}>\hat{l}=l_{\hat{k},\hat{n}}\geq\hat{k}$. Thus, it follows from the monotonicity of $\left\|x^\dagger- x_{k}^{(\hat{l})}\right\|^2$ with respect to $k$ that
\[
\left\|x^\dagger- x_{k_{\hat{l}}}^{(\hat{l})}\right\|^2\leq\left\|x^\dagger- x_{\hat{k}}^{(\hat{l})}\right\|^2<\varepsilon_0.
\]
 This leads to a contradiction to (\ref{equation 4-31}) with $l=\hat{l}$.

 Next, we turn to the discussion of $\{l_{k,n}\}$ and $(\hat{x}_k)$, for each $k=0,1,\cdots,$  which lead to
(i) and (ii) hold. In this case, we adopt a diagonal argument. For $k=0$, let $(\hat{x}_0)=(x_0)$ as well as $l_{0,n}=n$ for all $n$. Since $x_0^{(n)}=x_0$,
then (ii) holds automatically for $k=0$.

Assume that we have already constructed $\{l_{k,n}\}$  and $(\hat{x}_k)$  for all $0\leq k\leq m$. Then, we will focus on the definition of $\{l_{m+1,n}\}$ and $(\hat{x}_{m+1})$. According to the combination parameter $\lambda_m^{(l_{m,n})}$ and the integer $i_m^{(l_{m,n})}$ involved in the construction of $\left(x_{m+1}^{(l_{m,n})}\right)$. We get
\[
0\leq\lambda_m^{(l_{m,n})}\leq\frac{m}{m+\alpha}~~~~and~~~~0\leq i_m^{(l_{m,n})}\leq mj_{\max}~~\forall~n,
\]
where the constant $\alpha\geq3$. Choose a subsequence of $\{l_{m,n}\}$, denoted by $\{l_{m+1,n}\}$, which satisfies
\begin{equation}\label{equation 4-33}
\lim\limits_{n\rightarrow+\infty}\lambda_m^{(l_{m+1,n})}=\hat{\lambda}_m~~~and~~~i_m^{(l_{m+1,n})}=\hat{i}_m~~\forall~ n,
\end{equation}
for some number $0\leq\hat{\lambda}_m\leq m/(m+\alpha)$ and some integer $\hat{i}_m$. Define
\[
\hat{z}_m=\hat{x}_m+\hat{\lambda}_m(\hat{x}_m-\hat{x}_{m-1}).
\]
It follows from the induction hypothesis and (\ref{equation 4-33}) that
\begin{equation}\label{equation 4-34}
z_m^{(l_{m+1,n})}\rightarrow\hat{z}_m~~as~~n\rightarrow+\infty.
\end{equation}
Using (\ref{equation 4-33}), the continuity of $F, F'$, and the discussion in the proof of Lemma \ref{lemma 4-10}, there holds
\begin{equation}\label{equation 4-35}
t_{m,i}^{(l_{m+1,n})}\rightarrow\hat{t}_{m,i}~~as~~n\rightarrow+\infty,~~~\forall i\in I_m.
\end{equation}
Define
\[
\hat{x}_{m+1}=\hat{z}_{m}-\sum\limits_{i\in I_{m}}\hat{t}_{m,i} F'(\hat{z}_{m})^*(F(\hat{z}_{m})-y).
\]
According to (\ref{equation 4-34}) and (\ref{equation 4-35}), it follows that
\begin{equation}\label{equation 4-36}
\lim\limits_{n\rightarrow+\infty}x_{m+1}^{(l_{m+1,n})}=\hat{x}_{m+1}.
\end{equation}
Thus we complete the construction of $\{l_{m+1,n}\}$ as well as $(\hat{x}_{m+1})$.

We also need to prove that $\hat{\lambda}_m$ satisfies the requirements which in order to guarantee that the generated sequence is indeed in $\Gamma_{\eta, q}(x_0)$. It follows from the definition of the sequences in $\Gamma_{\eta, q}(x_0)$ that
\[
\lambda_{m}^{(l_{m+1,n})}(\lambda_{m}^{(l_{m+1,n})}+1)\left\|x_{m}^{(l_{m+1,n})}-x_{m-1}^{(l_{m+1,n})}\right\|^2
\leq\frac{\Psi^2}{\mu c_F^2}\|F(z_{m}^{(l_{m+1,n})})-y\|^2,
\]
as well as
\[
0\leq\lambda_{m}^{(l_{m+1,n})}\leq\min \left\{ {\frac{{q\left( i_m^{(l_{m+1,n})} \right)}}{{\left\| {x_m^{(l_{m+1,n})}  - x_{m -1}^{(l_{m+1,n})} } \right\|}},\frac{m}{m+\alpha}} \right\},
\]
with $1\leq i_m^{(l_{m+1,n})}-i_{m-1}^{(l_{m+1,n})}\leq j_{\max}$ for all $n$. Combining with (\ref{equation 4-33}), (\ref{equation 4-34}) and (\ref{equation 4-35}), as $n\rightarrow+\infty$  in the above two inequalities we can conclude that
\[
\hat{\lambda}_{m}(\hat{\lambda}_{m}+1)\left\|\hat{x}_{m}-\hat{x}_{m-1}\right\|^2
\leq\frac{\Psi^2}{\mu c_F^2}\|F(\hat{z}_{m})-y\|^2,
\]
and
\[
0\leq\hat{\lambda}_{m}\leq\min \left\{ {\frac{{q\left( \hat{i}_m \right)}}{{\left\| {\hat{x}_m  - \hat{x}_{m -1}} \right\|}},\frac{m}{m+\alpha}} \right\}.
\]
Assume that $\{l_{m,n}\}$ was chosen so that
$i_{m-1}^{(l_{m,n})}=\hat{i}_{m-1}$ for all $n$ and some integer $\hat{i}_{m-1}$. It follows that $\{l_{m+1,n}\}$  is a subsequence of $\{l_{m,n}\}$, then $1\leq\hat{i}_{m}-\hat{i}_{m-1}= i_m^{(l_{m+1,n})}-i_{m-1}^{(l_{m+1,n})}\leq j_{\max}$. According to $i_0^{(l)}=0$  for all $l$, we have $\hat{i}_0=0$. Therefore, the proof is complete.
\end{proof}
\end{prop}

\begin{lemma}\label{lemma 4-13}
Assume that $\{y^{\delta_n}\}$ be a sequence with noisy data satisfying $\|y^{\delta_n}-y\|\leq\delta_n$, where $\delta_n\rightarrow0$ as $n\rightarrow+\infty$. Let the combination parameters $\{\lambda_k^{\delta_n}\}$ are defined by DBTS method with $\{i_0^{\delta_n}\}=0$. Then, pick a fixed integer $k\geq0$, there exist a sequence $\{(x_i)\}\in \Gamma_{\eta, q}(x_0)$ such that
\[
x_i^{\delta_n}\rightarrow x_i,~~z_i^{\delta_n}\rightarrow z_i~~and ~~F(x_i^{\delta_n})\rightarrow F(x_i)~~as~n\rightarrow+\infty
\]
for all $0\leq i\leq k$.
\begin{proof}
 We prove the assertion by induction discussion on $k$. For $k=0$, we have $x_0^{\delta_n}=x_0$, $\lambda_0^{\delta_n}=0$ and then $F(x_0^{\delta_n})\rightarrow F(x_0)$, $z_0^{\delta_n}=z_0=x_0$. Now, assumed that the assertion is true for $k=m$, i.e., for all $0\leq i\leq m$, there exist a sequence $\{(x_i)\}\in \Gamma_{\eta, q}(x_0)$ which satisfies
\[
x_i^{\delta_n}\rightarrow x_i,~~z_i^{\delta_n}\rightarrow z_i~~and ~~F(x_i^{\delta_n})\rightarrow F(x_i)~~as~n\rightarrow+\infty.
\]

 Next, we show that the result with $i=m+1$ is also true. Without losing generality, we will get a sequence from $\Gamma_{\eta, q}(x_0)$  by maintaining the former $m+1$ terms in $\{(x_i)\}$ and modifying the rest terms. It is necessary to redefine $x_{m+1}$ because we
can apply the DBTS method with $\lambda_i=0$  for $i\leq m+1$ to generate the rest terms.

Furthermore, the combination parameter $\lambda_m^{\delta_n}$ generated by Algorithm \ref{algorithm 4-1} such that
\begin{equation}\label{equation 4-37}
\lambda_m^{\delta_n}(\lambda_m^{\delta_n}+1)\left\|x_m^{\delta_n}-x_{m-1}^{\delta_n}\right\|^2
\leq\frac{\Psi^2}{\mu c_F^2}\|F(z_m^{\delta_n})-y^{\delta_n}\|^2,
\end{equation}
as well as
\begin{equation}\label{equation 4-38}
\begin{array}{l}
0\leq\lambda_m^{\delta_n}\leq\min \left\{ {\frac{{q\left( i_m^{\delta_n} \right)}}{{\left\| {x_m^{\delta_n}  - x_{m-1}^{\delta_n}} \right\|}},\frac{m}{m+\alpha}} \right\}~or\\
 \lambda_m^{\delta_n}=\min\left\{\sqrt{\frac{(\Psi\tau\delta_n)^2}{\mu c_F^2\|{x_m^{\delta_n}-x_{m-1}^{\delta_n}}\|^2}+\frac{1}{4}}-\frac{1}{2},\frac{m}{m+\alpha}\right\}
\end{array}
\end{equation}
with $1\leq i_m^{\delta_{n}}-i_{m-1}^{\delta_{n}}\leq j_{\max}$. By defining a subsequence of $\{y^{\delta_n}\}$ if necessary, then we obtain
\begin{equation}\label{equation 4-39}
\lim\limits_{n\rightarrow+\infty}\lambda_m^{\delta_n}=\lambda_m~~and~~i_m^{\delta_n}=i_m~~for~ all~ n
\end{equation}
for some figure $0\leq\lambda_m\leq m/(m+\alpha)$ and some integer $i_m$. Next,  define $z_m,t_{m,i}$ with $i\in I_m$ as well as $x_{m+1}$ by (\ref{equation 3-3}) with $k=m$. It follows (\ref{equation 4-39}) and the proof of Lemma \ref{lemma 4-10}, we thus obtain
$$
z_m^{\delta_n}\rightarrow z_m~~and~~x_{m+1}^{\delta_n}\rightarrow x_{m+1}~~as~n\rightarrow+\infty,
$$
and
$$
F(x_{m+1}^{\delta_n})\rightarrow F(x_{m+1})~~as~n\rightarrow+\infty.
$$
Using the
induction hypothesis and by taking $n\rightarrow+\infty$ in (\ref{equation 4-37}), we can get that $\lambda_m$ and $i_m$ meet the requirements in the definition of $\Gamma_{\eta, q}(x_0)$. The proof is therefore complete.
\end{proof}
\end{lemma}

\begin{thm}\label{theorem 4-14}
Let the combination parameters $\{\lambda_k^{\delta_n}\}$ are defined by DBTS method with $i_0^{\delta_n}=0$. Find the integer $k_*$ which determined by the discrepancy principle (\ref{equation 3-10}). If  $\mathcal{N}\left(F'(x^\dagger)\right)\subset\mathcal{N}\left(F'(x)\right)$ for all $x\in B_{4\rho}(x^\dagger)$, then there holds
\begin{equation}\label{equation 4-40}
\lim\limits_{\delta\rightarrow0}\|x_{k_*}^\delta-x^\dagger\|=0.
\end{equation}
\begin{proof}
By contradiction, suppose that (\ref{equation 4-40}) is not true, then there exist a figure $\varepsilon>0$ as well as a subsequence $\{y^{\delta_n}\}$ of $\{y^{\delta}\}$ which satisfying $\|y^{\delta_n}-y\|\leq\delta_n$ with $\delta_n\rightarrow0$ as $n\rightarrow+\infty$ such that
\begin{equation}\label{equation 4-41}
\|x_{k_*}^{\delta_n}-x^\dagger\|\geq\varepsilon~~for~all~n.
\end{equation}
Hence, it remains to consider the following two cases.\\
\indent $Case~1.$ Assume that $\{k_{\delta_n}\}$  has a finite cluster point $\hat{k}$. By taking a subsequence if necessary and
using Lemma \ref{lemma 4-13}, we get $k_{\delta_n}=\hat{k}$ for all $n$ and there exist a sequence $\{(x_k)\}\in \Gamma_{\eta, q}(x_0)$  satisfies
\[
x_{\hat{k}}^{\delta_n}\rightarrow x_{\hat{k}},~~as~n\rightarrow+\infty.
\]
It can be seen via Theorem \ref{theorem 4-11} that we can obtain $x_{\hat{k}}=x^\dagger$. There holds
\[
x_{k_{*}}^{\delta_n}\rightarrow x^\dagger,~~as~n\rightarrow+\infty,
\]
i.e., $\|x_{k_*}^{\delta_n}-x^\dagger\|\rightarrow0$ as $n\rightarrow+\infty$. This leads to a contradiction.

 \indent $Case~2.$ Assume that $\lim\limits_{n\rightarrow+\infty}k_{\delta_n}=+\infty$. It follows Proposition \ref{proposition 4-12} that there exist an integer $k(\varepsilon)$ which satisfies
\begin{equation}\label{equation 4-42}
\|x_{k(\varepsilon)}-x^\dagger\|<\varepsilon~~for~any~\{(x_k)\}\in \Gamma_{\eta, q}(x_0).
\end{equation}
For this $k(\varepsilon)$, by taking a subsequence if necessary and using Lemma \ref{lemma 4-13}, we can pick $\{(x_k)\}\in \Gamma_{\eta, q}(x_0)$ such that
\begin{equation}\label{equation 4-43}
x_{k}^{\delta_n}\rightarrow x_k~~as~n\rightarrow+\infty,~~\forall~0\leq k\leq k(\varepsilon).
\end{equation}
It follows from $\lim\limits_{n\rightarrow+\infty}k_{*}=+\infty$, that we have $k_{*}>k(\varepsilon)$ for $n$ which is large enough. Then, according to (\ref{equation 4-10})  in Proposition \ref{proposition 4-2} we get
\[
\left\|x_{k_*}^{\delta_n}-x^\dagger\right\|\leq\left\|x_{k(\varepsilon)}^{\delta_n}-x^\dagger\right\|.
\]
Combining (\ref{equation 4-42}) and (\ref{equation 4-43}) as $n\rightarrow+\infty$, we arrive at
\[
\lim\limits_{n\rightarrow+\infty}\|x_{k(\varepsilon)}^{\delta_n}-x^\dagger\|=\|x_{k(\varepsilon)}-x^\dagger\|<\varepsilon.
\]
This is a contradiction to (\ref{equation 4-41}).
According to the above discussions, we therefore get (\ref{equation 4-40}) which yields the assertion.
\end{proof}
\end{thm}

\begin{remark}\label{remark 4-4}
In the proof of Theorem \ref{theorem 4-8}, Lemma \ref{lemma 4-10} Theorem \ref{theorem 4-11}, we take advantage of $\lambda_k^\delta$ with $\delta\rightarrow0$ by defining the hypothesis of uniqueness.
However, given the fact that it is not unique in the DBTS method, we
discuss the generated different cluster points above. In Theorem \ref{theorem 4-14} we removed this hypothesis by using a uniform convergence result established in Proposition \ref{proposition 4-12}.
\end{remark}

\section{Numerical simulations}\label{section-5}

\noindent In this section we carry out some numerical experiments with one-dimensional and two-dimensional inverse potential problem to test the good performance of the proposed TGSS method for solving the ill-posed system (\ref{equation 1-1}). Our simulations were done by using MATLAB R2013a on a LG computer with Intel Core i5-6500 CPU 3.20 GHz and 8.00 GB memory.

To comprehensively demonstrate the acceleration performance of TGSS method with two search directions,
we compare the proposed TGSS method with the algorithms associated with it. First, we explain the following abbreviations.

\begin{description}
  \item[1.] Land: The classical Landweber method (\ref{equation 1-3}) which
  we choice the  stepsize with a constant $\alpha_k^\delta=1$.
  \item[2.] TPG-DBTS: Two-point gradient method, whose iteration scheme is given by (\ref{equation 1-5}), and the parameters $\lambda_k^\delta$ are selected in DBTS method.
  \item[3.] TPG-Nes: Two-point gradient method with classic Nesterov acceleration scheme, as equation (\ref{equation 1-4}).
  \item[4.] SESOP:  Sequential subspace optimization method for the iterative solution, i.e., equation (\ref{equation 1-6}).
  \item[5.] TGSS-DBTS: Accelerated two-point gradient method based on the sequential subspace optimization with $\lambda_k^\delta$ are chosen in DBTS method.
  \item[6.] TGSS-Nes: Accelerated two-point gradient method based on the sequential subspace optimization with classic Nesterov acceleration scheme.
\end{description}

In order to better illustration, we define the following quantities for quantitatively analyze the reconstruction performance.
\begin{enumerate}
  \item Difference:
  \[
  x-x^\dagger,
  \]
  which is the difference between exact and the reconstructed solution.
  \item  Relative error (RE):
  \[
  RE=\frac{\|x-x^\dagger\|}{\|x^\dagger\|}.
  \]
  \item The rate of iterations:
  \[
  Rate(k_*)=\frac{k_*(\cdot)}{k_*(Land)},
  \]
  which illustrates the degree of acceleration on the aspect of iterations. Here, $k_*(\cdot)$ represents the number of iteration of corresponding  algorithms under the discrepancy principle, and
  $k_*(Land)$  is the number of iteration of Landweber  iteration under the same conditions.
  \item The rate of computation time:
  \[
  Rate(t)=\frac{t(\cdot)}{t(Land)},
  \]
  which illustrates the degree of acceleration on computation time. Here, $t(\cdot)$ represents the computation time of corresponding  algorithms, and
  $t(Land)$  is the computation time of Landweber iteration.
\end{enumerate}

 \noindent Next, we consider the following elliptic equation
\begin{equation}\label{equation 5-1}
\left\{
  \begin{array}{ll}
    -\Delta u+cu=f, & \hbox{in $\Omega$,} \\
    \frac{\partial u}{\partial n}=0, & \hbox{on $\Gamma$,}
  \end{array}
\right.
\end{equation}
from the measurement of $u$, where $\Omega\subset\mathbb{R}^d(d=1,2)$ is an open bounded domain with a Lipschitz boundary $\Gamma$ as well as $f\in L^2(\Omega)$. We assume that the sought solution $c\in L^2(\Omega)$. Then the inverse problem can be equivalently described as solving a nonlinear operator equation
\[
F(c)=u,
\]
we define the nonlinear operator $F:\mathcal{D}(F)\mapsto L^2(\Omega)$. The domain $\mathcal{D}(F)$ for
\[
\mathcal{D}(F):=\left\{c\in L^2(\Omega):\|c-\hat{c}\|_{L^2(\Omega)}\leq \zeta_0~for~ some~\hat{c}\geq0, a.e.\right\}
\]
as the admissible set of $F$. As we know that $F$ is well-defined for some positive constant $\zeta_0$ as well as Fr\'{e}chet differentiable \cite{26}. The Fr\'{e}chet derivative of $F$ and its adjoint can be calculated as follows
\[
F'(c)q=-A(c)^{-1}(qF(c))~=and~F'(c)^*\omega=-u(c)A(c)^{-1}\omega,
\]
for $q,\omega\in L^2(\Omega)$. Here $A(c):H^2(\Omega)\cap H_0^1(\Omega)\rightarrow L^2(\Omega)$ is defined by $A(c)u=-\Delta u+cu$ \cite{27}.

 By using $u^\dagger=F(c^\dagger)$ with the true parameter $c^\dagger$, we get the exact data $u^\dagger$. We consider Gaussian noisy data models where we simulate noisy data $u^\delta$ in the following way
\[
u^\delta=u^\dagger+\delta\cdot n,
\]
where $\delta$ represents the noise level and $n$ represents the random variable satisfying the standard normal distribution. In this paper, the forward
operator of the elliptic equation is discretized using finite elements on a uniform grid which we refer to \cite{28}. The implementation of corresponding experiment with reference to the MATLAB package shared by Bangti Jin (http://www.uni-graz.at/$\backsim$ clason/codes/l1fitting.zip).

\subsection{1-D inverse potential problem}\label{subsection 5-1}
\noindent We first consider the one-dimensional case of above problem with the following information.

\begin{itemize}
  \item Let $\Omega=[-1,1],~f(x)\equiv1$ as well as true parameter
  \[
  c^\dagger(x)=1-\cos(\pi x)
  \]
  for equation (\ref{equation 5-1}).
  \item  The exact data $u$ are obtained by solving the forward model using the standard piecewise linear finite element method with mesh $N=256$ and mesh size $h=1/N$.
  \item For comparison purposes, we choose three noise levels, $\delta=0.1\%$, $\delta=0.01\%$ and $\delta=0.001\%$.
  \item Assume that $\eta=0.1$. Since $\tau$ needs to satisfy $\tau>\frac{1+\eta}{1-\eta}$ in SESOP method and TGSS method, as well as the condition $\tau>2\frac{1+\eta}{1-2\eta}$ in TPG method which can guarantee the corresponding convergence. Thus, we take $\tau=2.8$,  $\mu=1.01$ and $c_F=0.1$.
  \item As for the selection of combination parameters $\lambda_k^\delta$ in DBTS method set $i_0=2$, $j_{\max}=1$ and $q(i)=4/i^{1.1}$. In Nesterov  acceleration scheme we take $\lambda_k^\delta=\frac{k-1}{k+\alpha-1}$ with $\alpha=3$.
\end{itemize}

\begin{figure}[htb]
\centering
{
\subfigure[]{
\label{figure 1-1}
\includegraphics[scale=0.34]{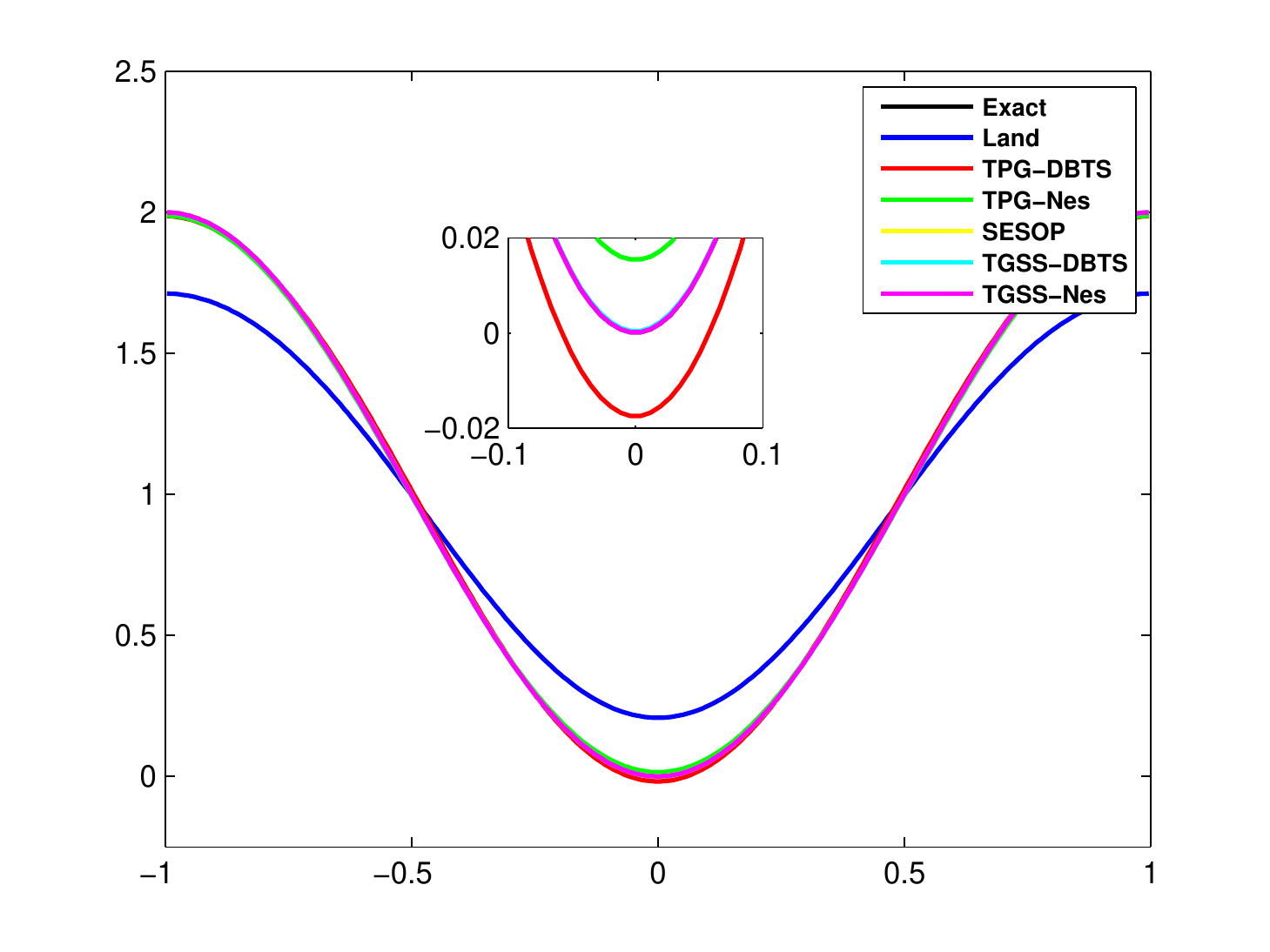} }
\subfigure[]{
\label{figure 1-2}
\includegraphics[scale=0.34]{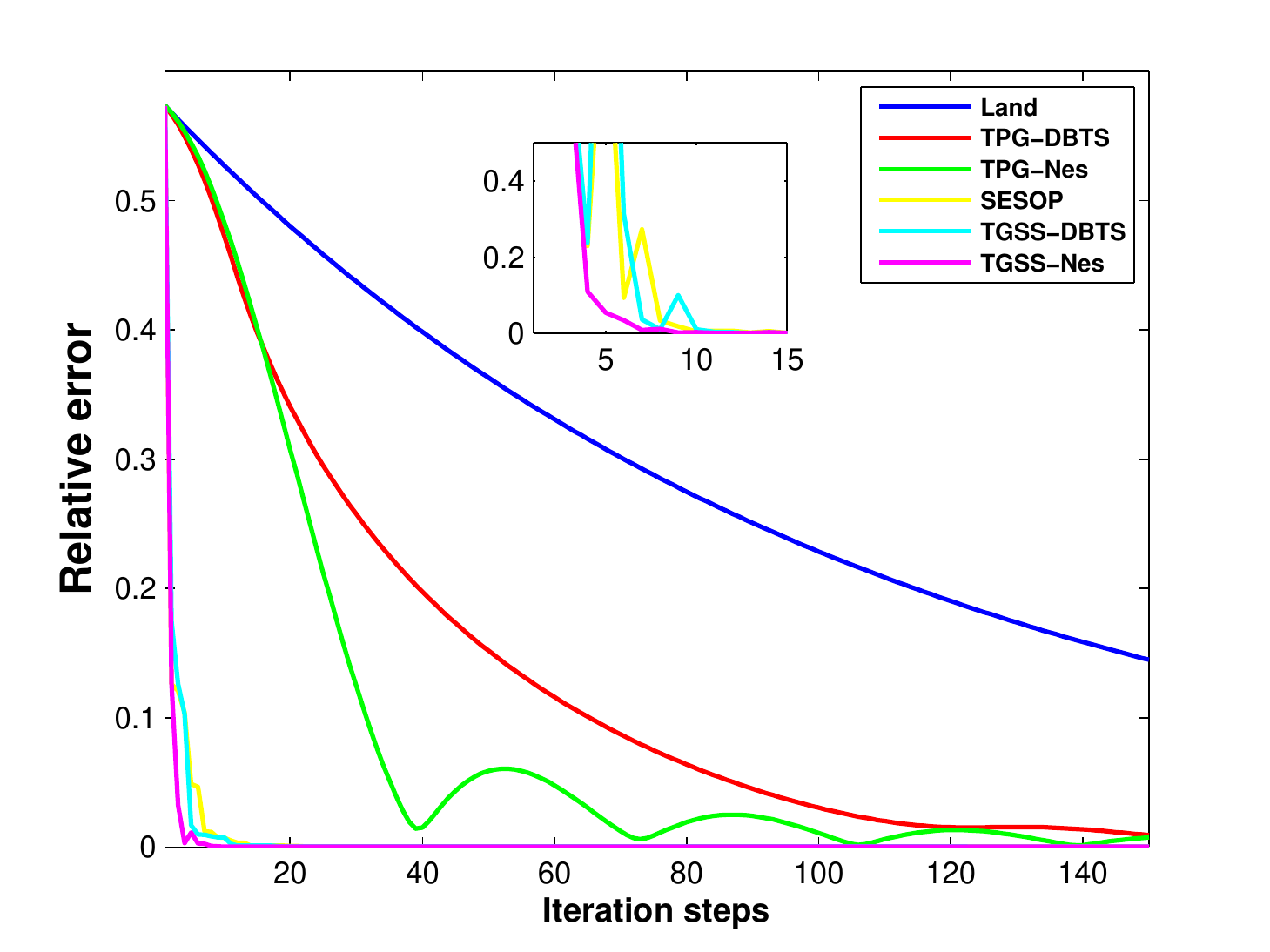} }
\subfigure[]{
\label{figure 1-3}
\includegraphics[scale=0.34]{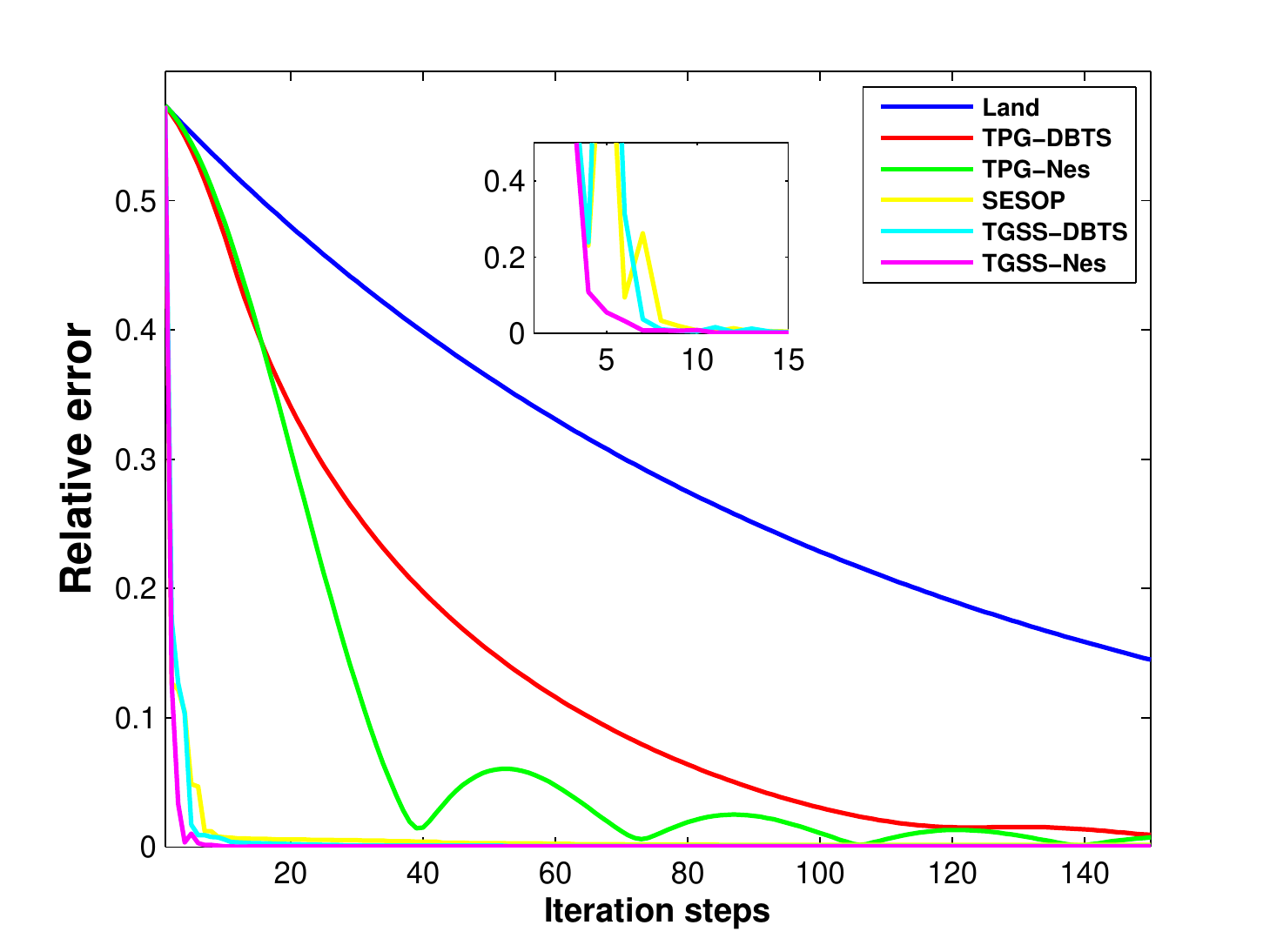} }
}
\caption{Reconstructions of each method for 1-D problem to show reconstructed solution and relative errors in the first 150 iterations. (a) Reconstructed solution. (b) Evolution of the relative error with exact data. (c) Evolution of the relative error with noise level $\delta=0.1\%$.
}
\label{1D}
\end{figure}

 We first illustrate the phenomenon generated by the numerical experiment after 150 steps of iteration which can be clearly seen from figure \ref{1D}. As shown in figure \ref{figure 1-1}, we observe that the reconstructed solutions generated by algorithm TGSS are closer to the true solution than other algorithms.
Meanwhile, when referring to relative error (RE) curves, it can be seen the advantage of our algorithm in convergence behavior in this respect from the figure \ref{figure 1-2} and figure \ref{figure 1-3} which respectively represent the case with exact data and the case with noise level of $\delta=0.1\%$. They clearly show that the TGSS method makes the reconstruction error decrease dramatically.


\begin{table}
\footnotesize
\centering
\caption{
Comparisons between Land, TPG, SESOP and TGSS methods for 1-D inverse potential problem with noisy data.
}
\begin{threeparttable}
\begin{tabular}{l l l l l l l l l l l l l}
\hline
$\delta$& Methods& $k_*$& Time (s)&  RE & Rate($k_*$) & & Rate(t)&\\
\hline
\multirow{1}{*}{0.1$\%$} & $\textrm{\textrm{Land}}$& 2169 & 2.003 & ${5.25 \times 10^{-3}}$ & 100$\%$ &  & 100$\%$ & \\
&$\textrm{\textrm{TPG-DBTS}}$& 182 & 0.180 & ${2.33 \times 10^{-3}}$ & 8.39$\%$ & & 8.99$\%$ & \\
&$\textrm{\textrm{TPG-Nes}}$& 105 & 0.101 & ${2.46 \times 10^{-3}}$ & 4.84$\%$ & & 5.04$\%$ & \\
&$\textrm{\textrm{SESOP}}$& 37 & 0.143 & ${4.24 \times 10^{-3}}$ & 1.71$\%$ & & 7.14$\%$ & \\
&$\textrm{\textrm{TGSS-DBTS}}$& 14 & 0.039 & ${3.04 \times 10^{-3}}$ & 0.65$\%$ & & 1.95$\%$ & \\
&$\textrm{\textrm{TGSS-Nes}}$& 10 & 0.025 & ${5.38 \times 10^{-4}}$ & 0.46$\%$ & & 1.25$\%$ & \\
\hline
\multirow{1}{*}{0.01$\%$} & $\textrm{\textrm{Land}}$& 5756 & 5.349 & ${6.97 \times 10^{-4}}$ & 100$\%$ & & 100$\%$ & \\
&$\textrm{\textrm{TPG-DBTS}}$& 233 & 0.231 & ${7.66 \times 10^{-4}}$ & 4.05$\%$ & & 4.32$\%$ & \\
&$\textrm{\textrm{TPG-Nes}}$& 336 & 0.316 & ${5.16 \times 10^{-4}}$ & 5.84$\%$ & & 5.91$\%$ & \\
&$\textrm{\textrm{SESOP}}$& 87 & 0.370 & ${5.33\times 10^{-4}}$ & 1.51$\%$ & & 6.92$\%$ & \\
&$\textrm{\textrm{TGSS-DBTS}}$& 24 & 0.068 & ${4.03 \times 10^{-4}}$ & 0.42$\%$ & & 1.27$\%$ & \\
&$\textrm{\textrm{TGSS-Nes}}$& 18 & 0.051 & ${1.34 \times 10^{-4}}$ & 0.31$\%$ & & 0.95$\%$ & \\
\hline
\multirow{1}{*}{0.001$\%$} & $\textrm{\textrm{Land}}$& 16960 & 15.763& $ {1.17 \times 10^{-4}}$ & 100$\%$ & & 100$\%$ & \\
&$\textrm{\textrm{TPG-DBTS}}$& 624 & 0.611 & ${5.35\times 10^{-5}}$ & 3.68$\%$ & & 3.88$\%$ & \\
&$\textrm{\textrm{TPG-Nes}}$& 732 & 0.685 & ${5.44 \times 10^{-5}}$ & 4.32$\%$ & & 4.35$\%$ & \\
&$\textrm{\textrm{SESOP}}$& 342 & 1.134 & ${1.12 \times 10^{-4}}$ & 2.02$\%$ & & 7.19$\%$ & \\
&$\textrm{\textrm{TGSS-DBTS}}$& 122 & 0.276 & ${1.08 \times 10^{-4}}$ & 0.72$\%$ & & 1.75$\%$ & \\
&$\textrm{\textrm{TGSS-Nes}}$& 103 & 0.177 & ${1.14 \times 10^{-5}}$ & 0.61$\%$ & & 1.12$\%$ & \\
\hline
\end{tabular}
\end{threeparttable}
\label{table 1-1}
\end{table}

In addition, three various noise levels are added to the generated exact data, respectively, to test the robustness of the proposed method (TGSS) to noise. Table \ref{table 1-1} summarizes detailed simulation results for each noise level. Under the same discrepancy principle, we can reach the following conclusions:
\begin{itemize}
  \item At all noise levels, compared with Land, TPG and SESOP methods, our proposed TGSS method leads to a strongly decrease of the iteration numbers and the overall computational time can be significantly reduced.
  \item We observe that TPG method has a good acceleration effects relative to Land, but when it combined with SESOP method (TGSS), the acceleration performance becomes more obvious.
  \item By using comparisons under the same conditions, i.e., the parameters $\lambda_k^\delta$ are chosen in the same way such as DBTS method or Nesterov acceleration scheme, we can get that TGSS method achieves more satisfactory acceleration effects than TPG method.
  \item TGSS-DBTS as well as TGSS-Nes methods seem to have similar convergence behavior which illustrate the validity of selecting $\lambda_k^\delta$ in DBTS method.
\end{itemize}

\subsection{2-D inverse potential problem}\label{subsection 5-2}
\noindent We now consider a 2-D case to illustrate the performance of  proposed methods on solving  inverse potential problem under different noise levels.

\begin{itemize}
  \item Let $\Omega=[-1,1]^2,~f(x_1,x_2)\equiv1$ as well as true parameter
  $$
  c^\dagger(x_1,x_2)=1+\cos(\pi x_1)\cos(\pi x_2)\chi_{\{|(x_1,x_2)|_\infty<1/2\}}
  $$
  for equation (\ref{equation 5-1}).
  \item  Use a sample with the two-dimensional standard piecewise linear finite element method with mesh $N=64\times64$ and mesh size $h=1/N^2$.
  \item Choose three noise levels, which are $\delta=2\%$, $\delta=1\%$ and $\delta=0.5\%$. Moreover, take $\eta=0.1$, $\tau=2.8$, $\mu=1.01$ and $c_F=0.1$.
  \item As for the selection of combination parameters $\lambda_k^\delta$ in DBTS method we set $i_0=2$, $j_{\max}=1$ and $q(i)=9/i^{1.1}$. In Nesterov acceleration scheme we take $\lambda_k^\delta=\frac{k-1}{k+\alpha-1}$ with $\alpha=9$.
\end{itemize}

\begin{figure}[htb]
\centering
{
\subfigure[]{
\label{figure 2-1}
\includegraphics[scale=0.34]{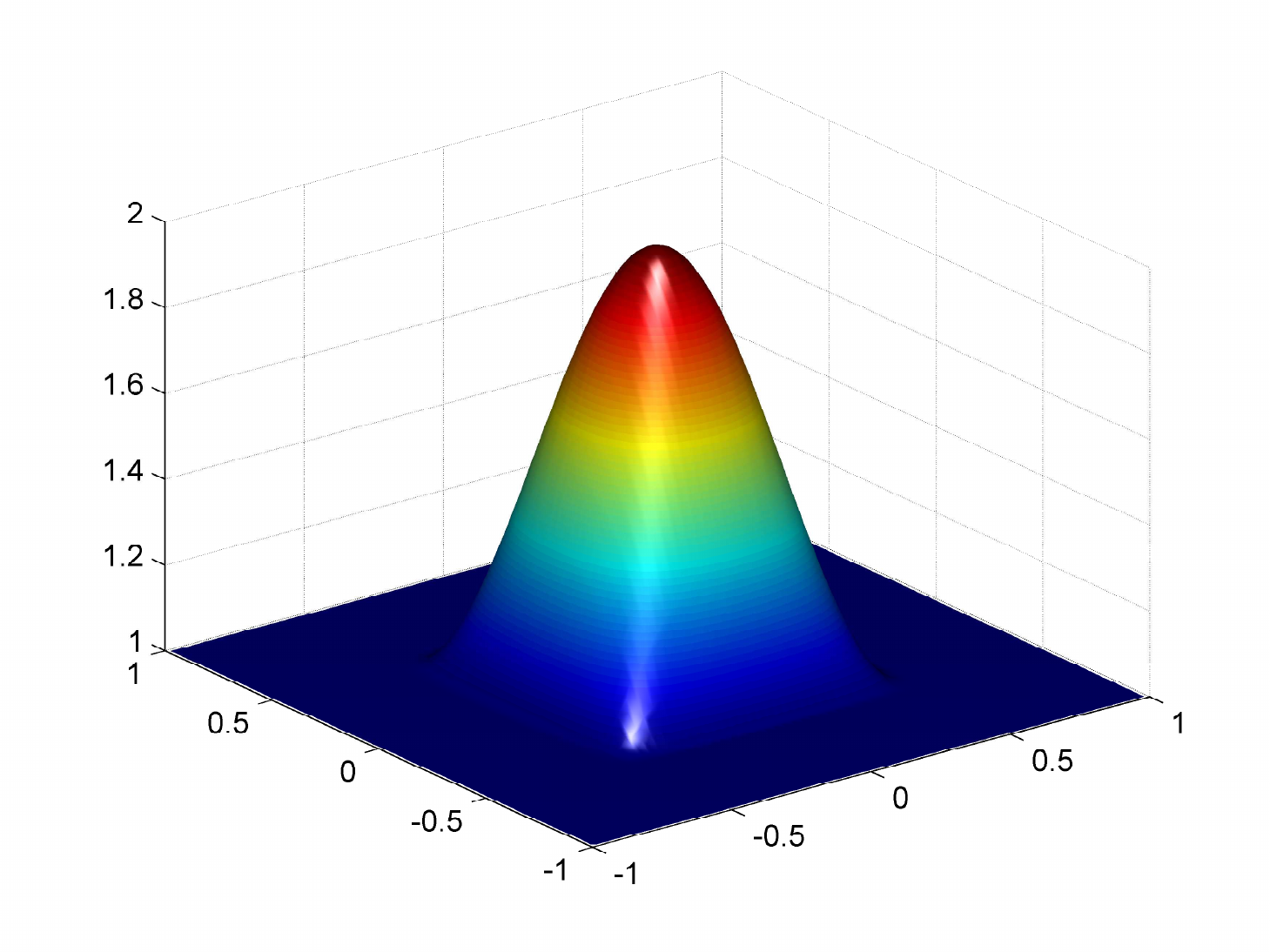} }
\subfigure[]{
\label{figure 2-2}
\includegraphics[scale=0.34]{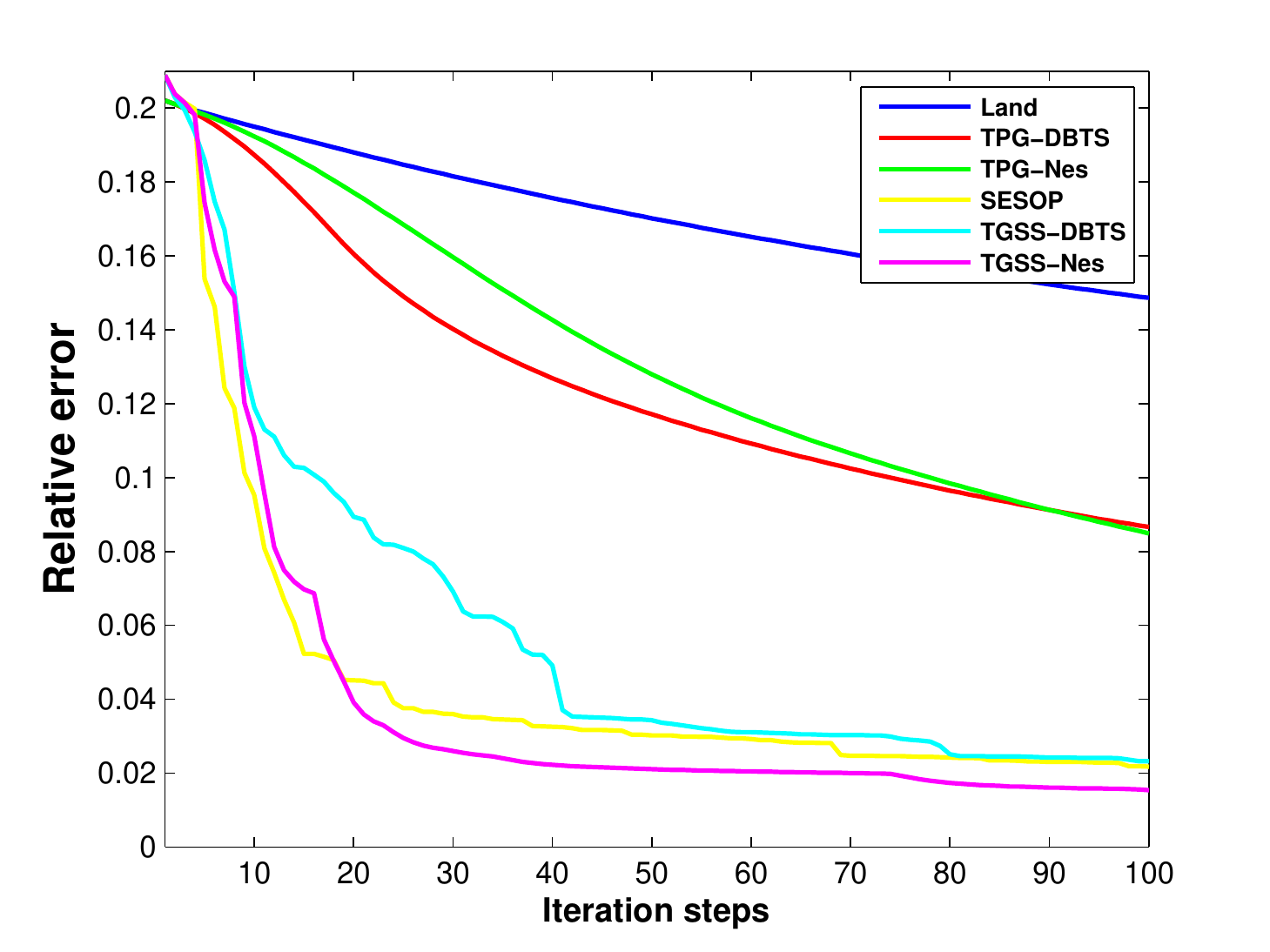} }
\subfigure[]{
\label{figure 2-3}
\includegraphics[scale=0.34]{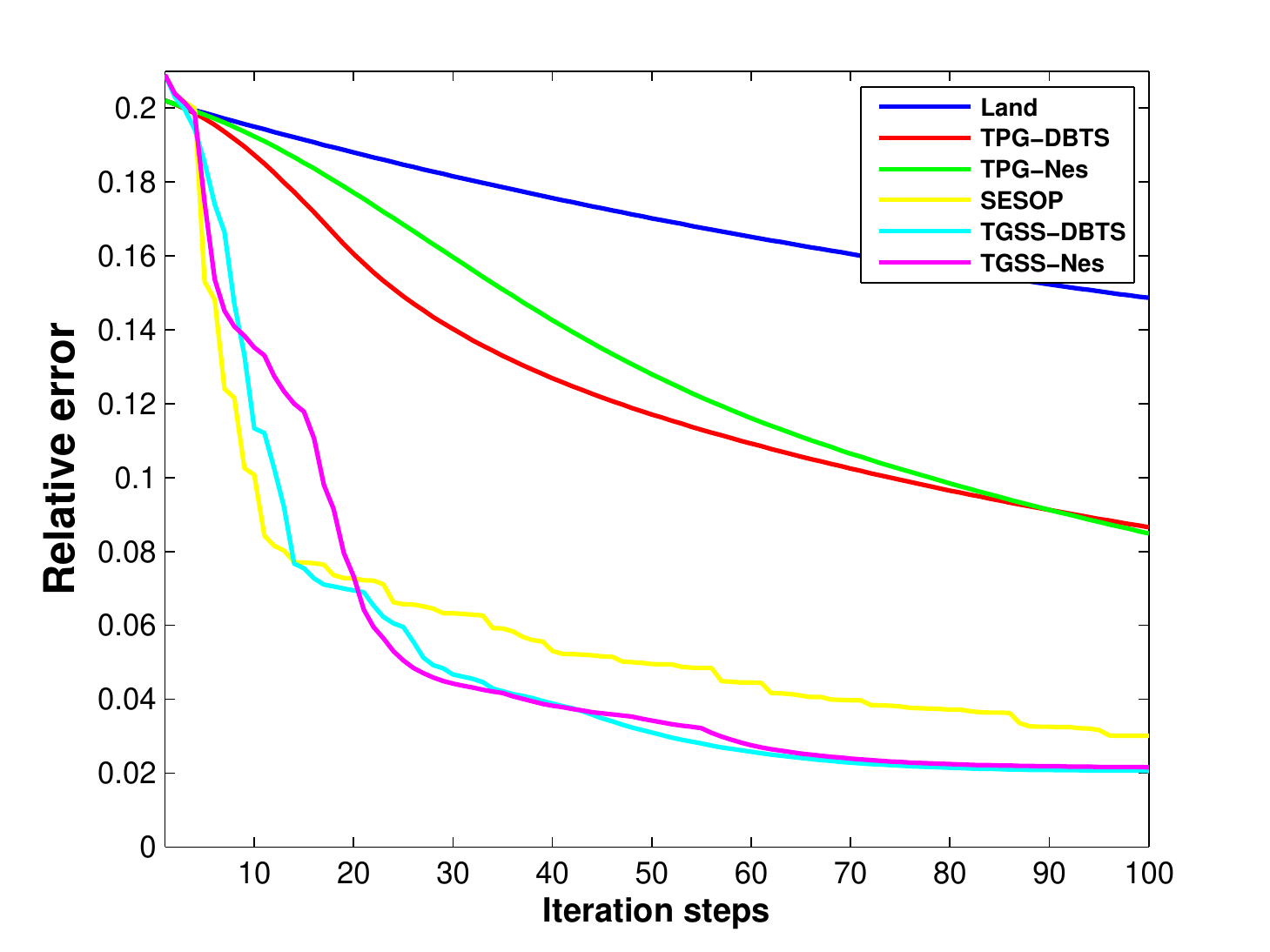} }
}
\caption{Experiment for 2-D problem to show true solution and the relative error evolution of each method in the first 100 iterations. (a) True solution. (b) Evolution of the relative error with exact data. (c) Evolution of the relative error with noise level $\delta=2\%$.
}
\label{2D-1}
\end{figure}

\begin{figure}[htb]
\centering
{
\subfigure[Reconstruction of Land]{
\label{figure 3-1}
\includegraphics[scale=0.34]{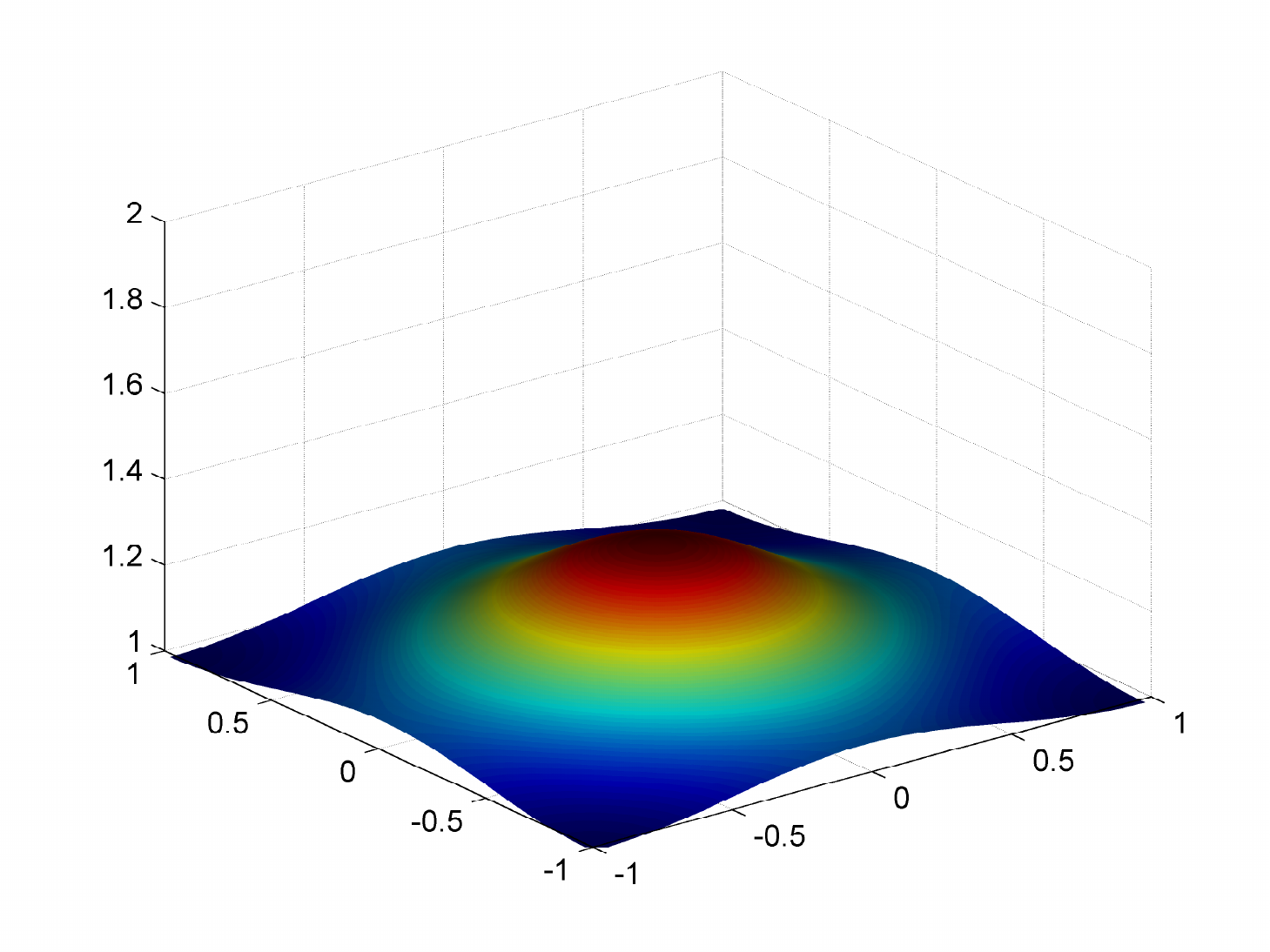} }
\subfigure[Reconstruction of TPG-DBTS]{
\label{figure 3-2}
\includegraphics[scale=0.34]{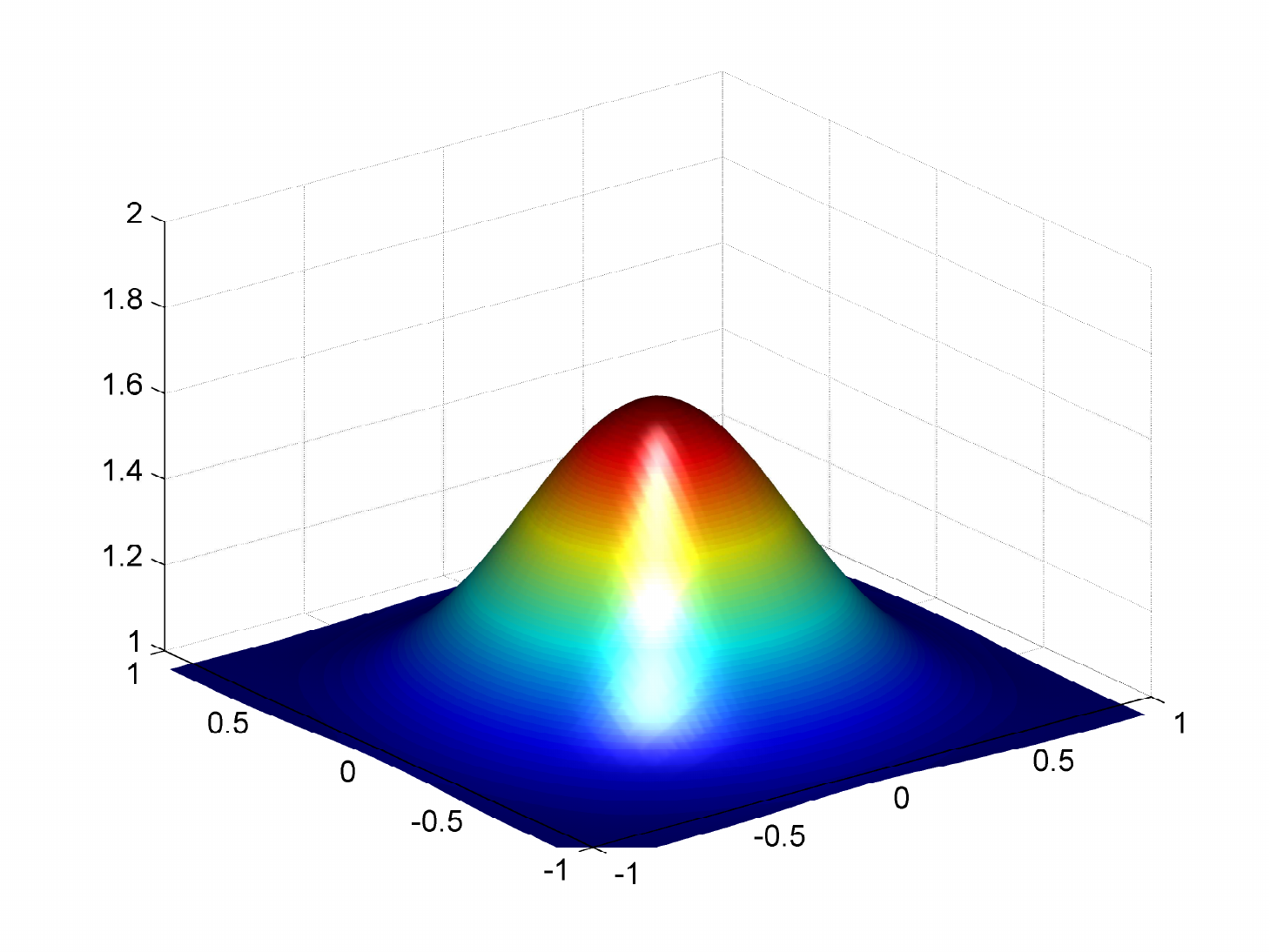} }
\subfigure[Reconstruction of TPG-Nes]{
\label{figure 3-3}
\includegraphics[scale=0.34]{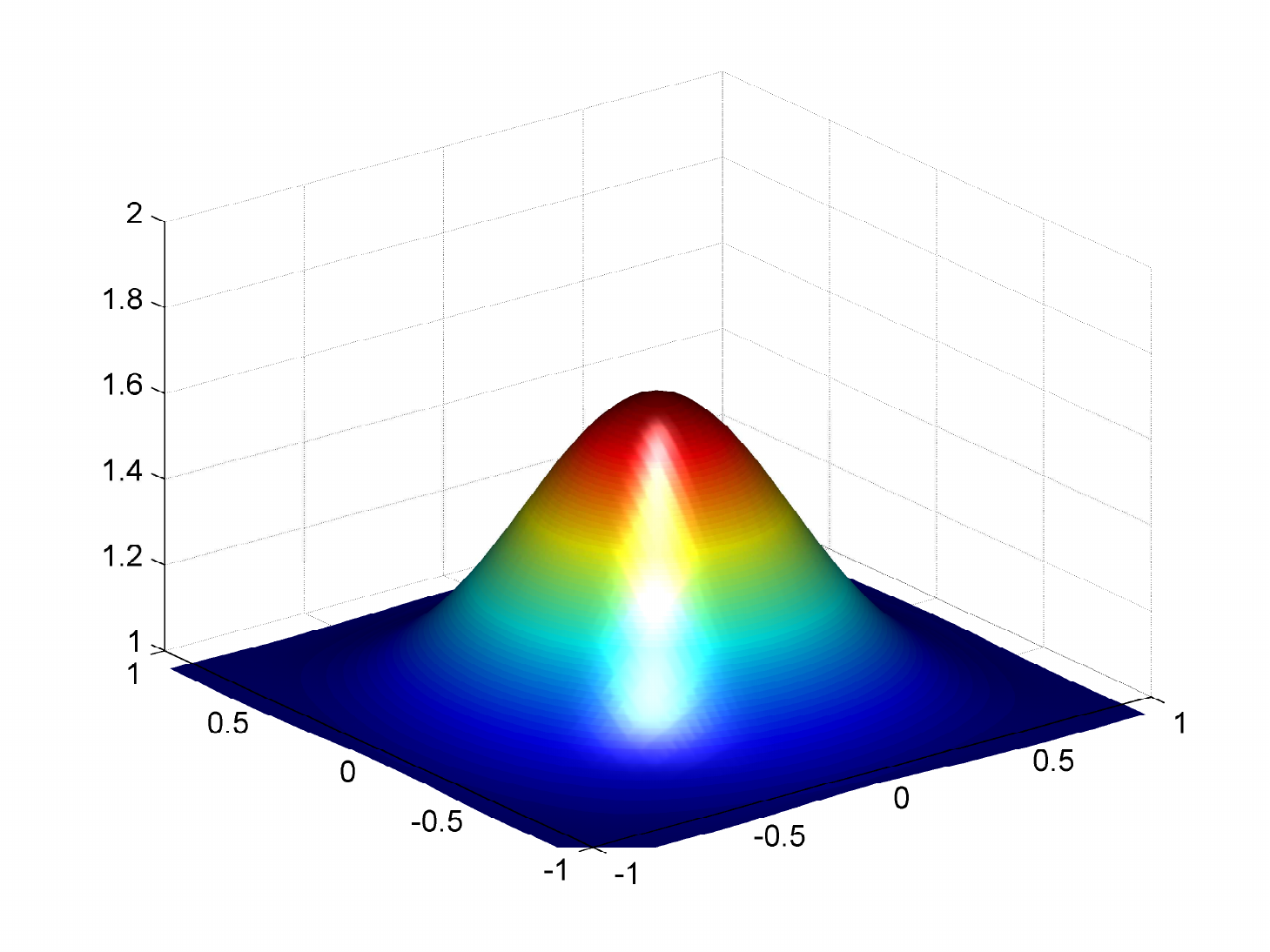} }
\subfigure[Reconstruction of SESOP]{
\label{figure 3-4}
\includegraphics[scale=0.34]{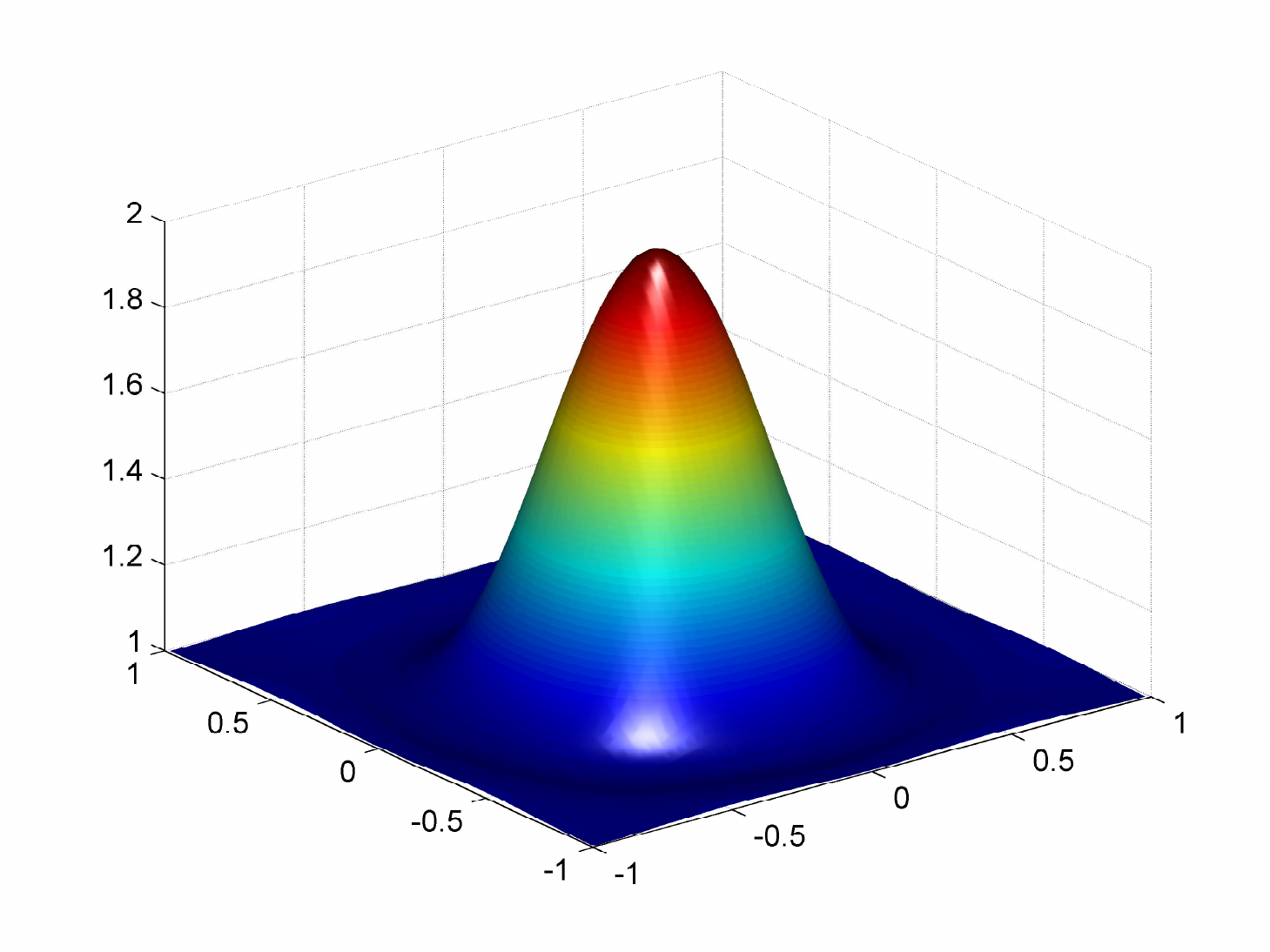} }
\subfigure[Reconstruction of TGSS-DBTS]{
\label{figure 3-5}
\includegraphics[scale=0.34]{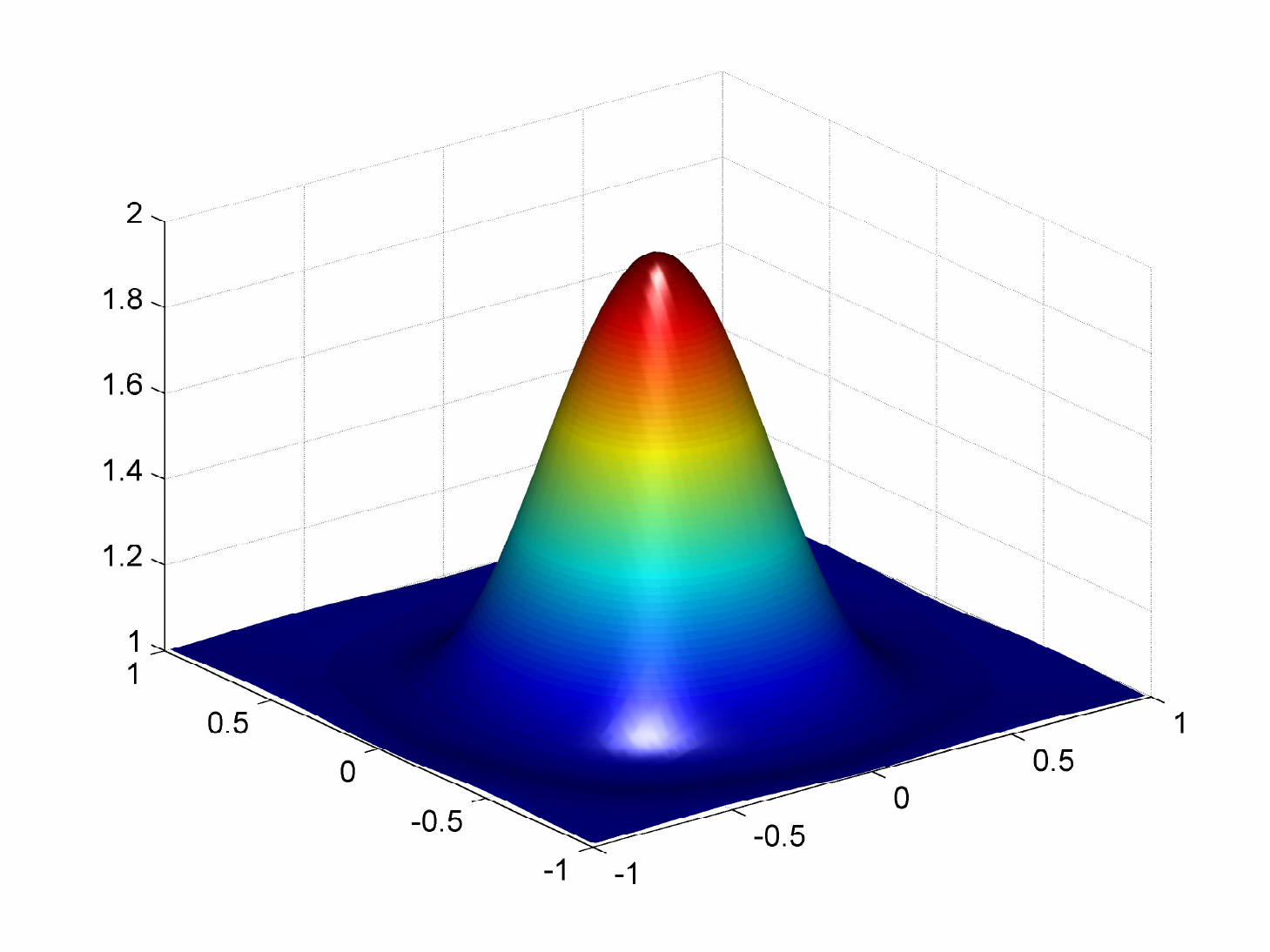} }
\subfigure[Reconstruction of TGSS-Nes]{
\label{figure 3-6}
\includegraphics[scale=0.34]{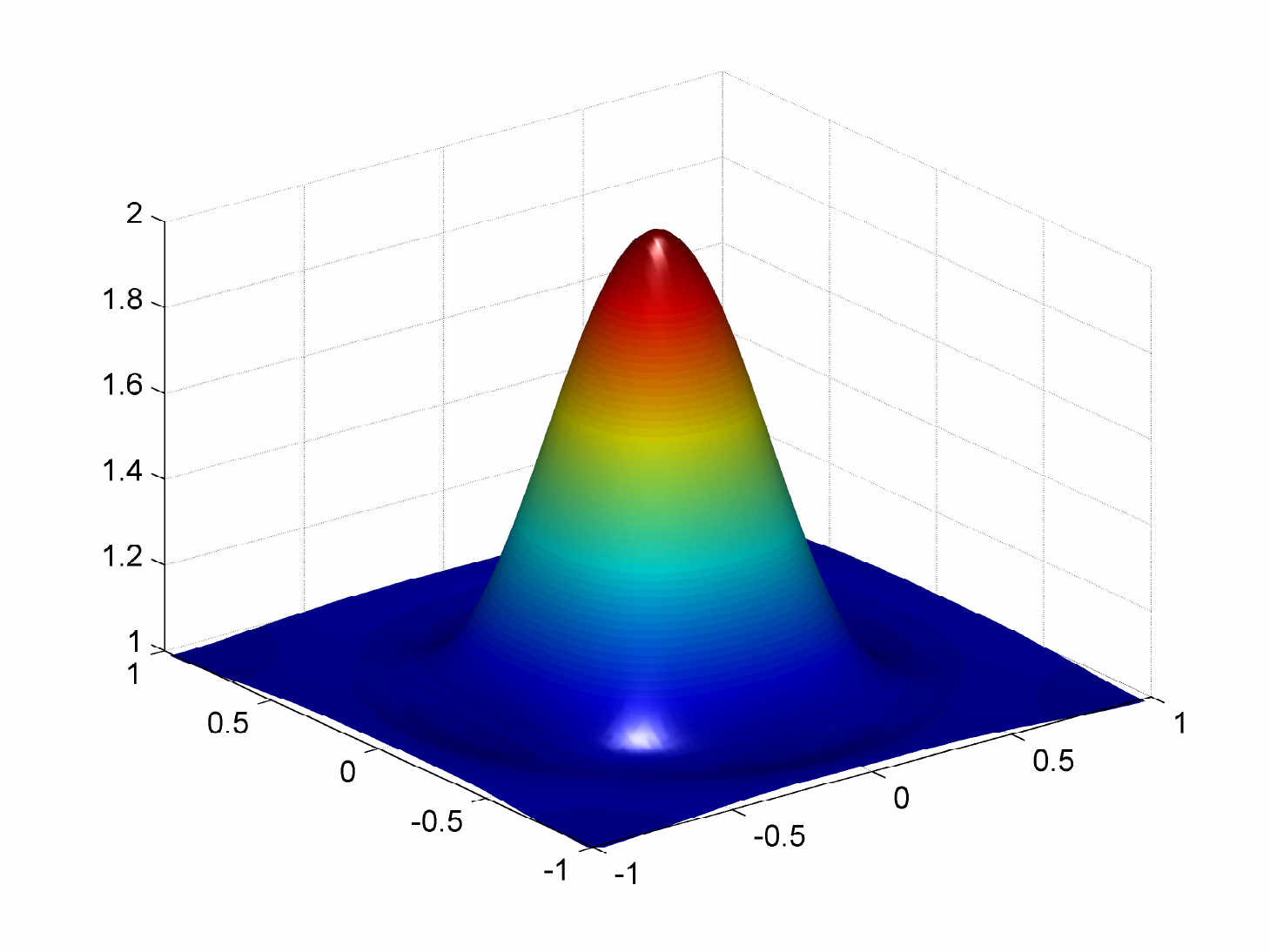} }
}
\caption{Experiment with exact data for 2-D  problem to show reconstructed solution of each method at 100th iteration. (a)-(f) Reconstructed solution of each method.
}
\label{2D-2}
\end{figure}

Compared with other methods, the 2-D numerical experiments depicted in Figure \ref{2D-1} and Figure \ref{2D-2} indicate that, in the exact data case the proposed method (TGSS) has better reconstruction effects. It follows Figure \ref{figure 2-1} (true solution) and Figure \ref{2D-2} (reconstructed solution of each method) that after 100 steps of iteration that our proposed method can obtain more accurate reconstruction results.  As can be
seen from Figure \ref{figure 2-2}, compared with other algorithms especially Land and TPG, the value of RE in our proposed method is going down faster.
Thus, we can conclude that TGSS method has a satisfactory convergence behavior in the case of exact data.

\begin{figure}[htb]
\centering
{
\subfigure[Difference by Land]{
\label{figure 4-1}
\includegraphics[scale=0.34]{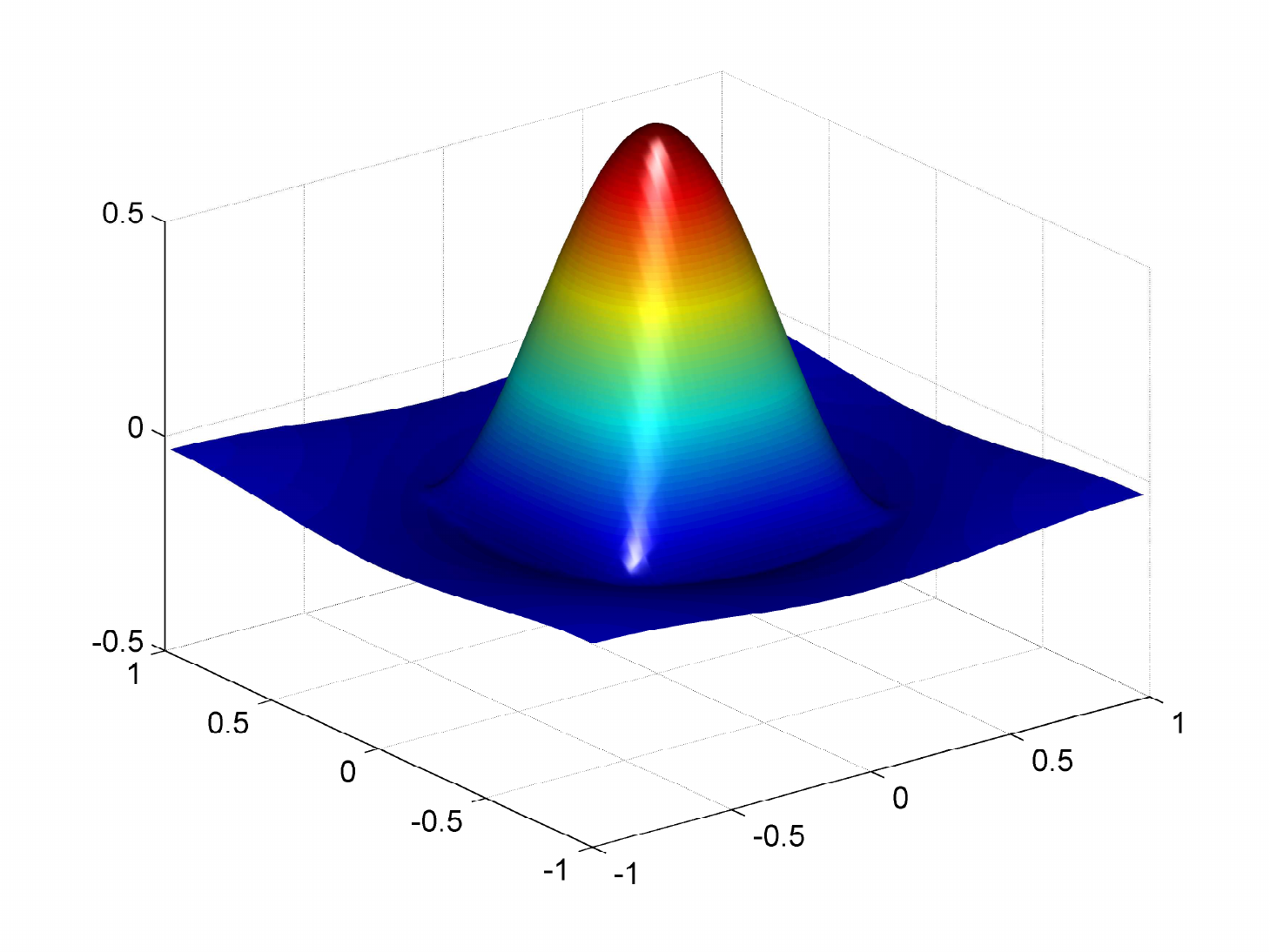} }
\subfigure[Difference by TPG-DBTS]{
\label{figure 4-2}
\includegraphics[scale=0.34]{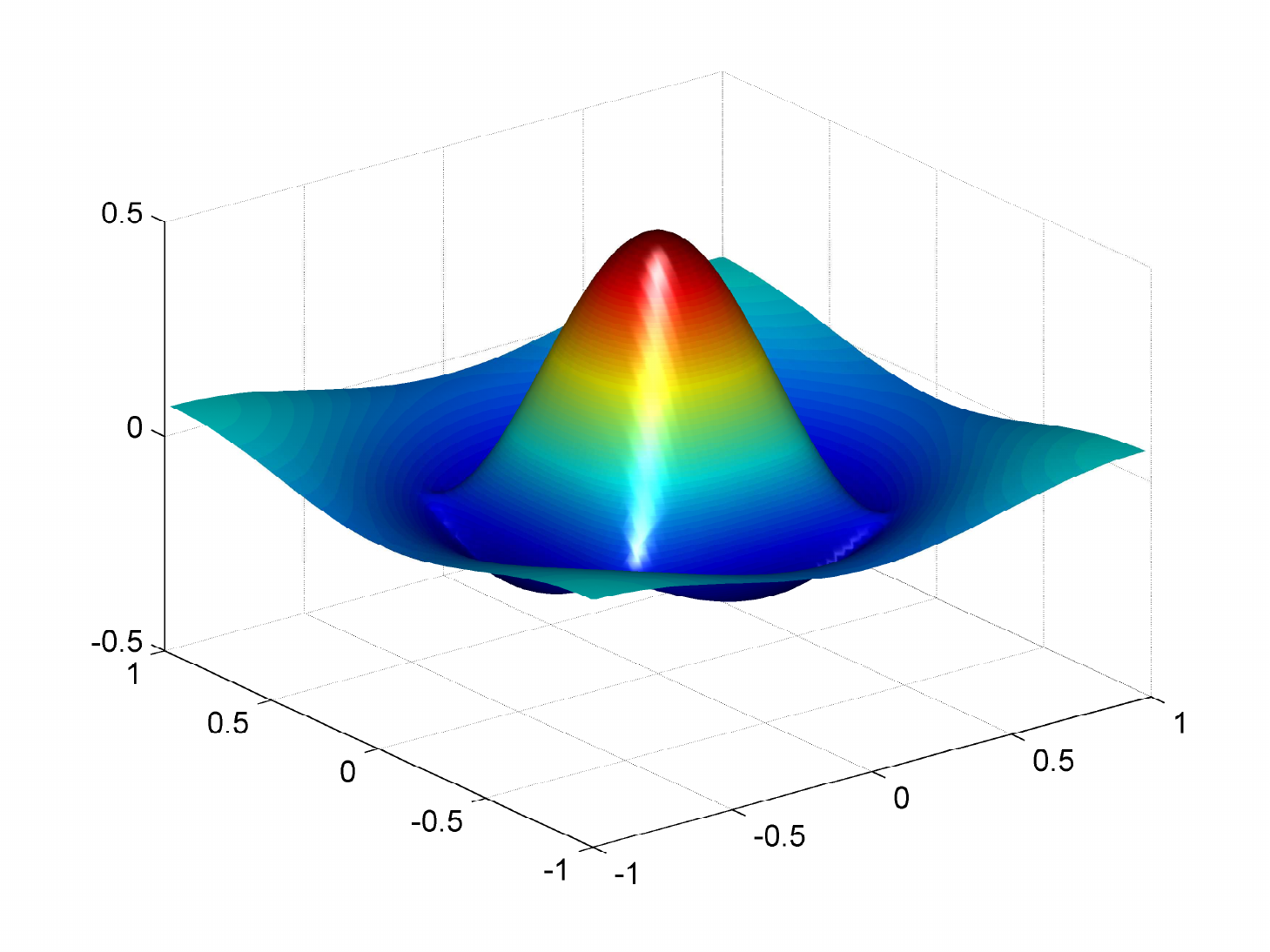} }
\subfigure[Difference by TPG-Nes]{
\label{figure 4-3}
\includegraphics[scale=0.34]{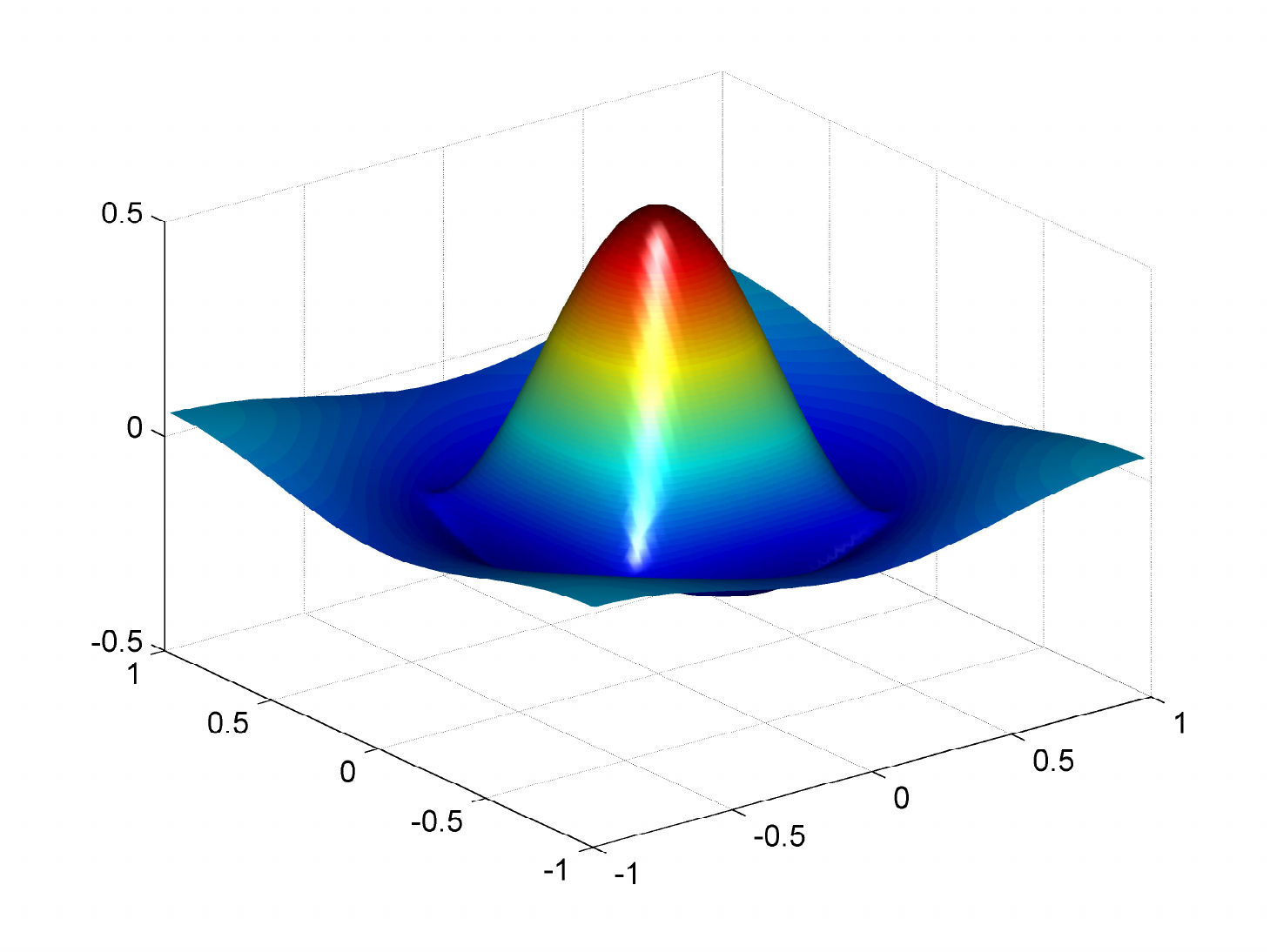} }
\subfigure[Difference by SESOP]{
\label{figure 4-4}
\includegraphics[scale=0.34]{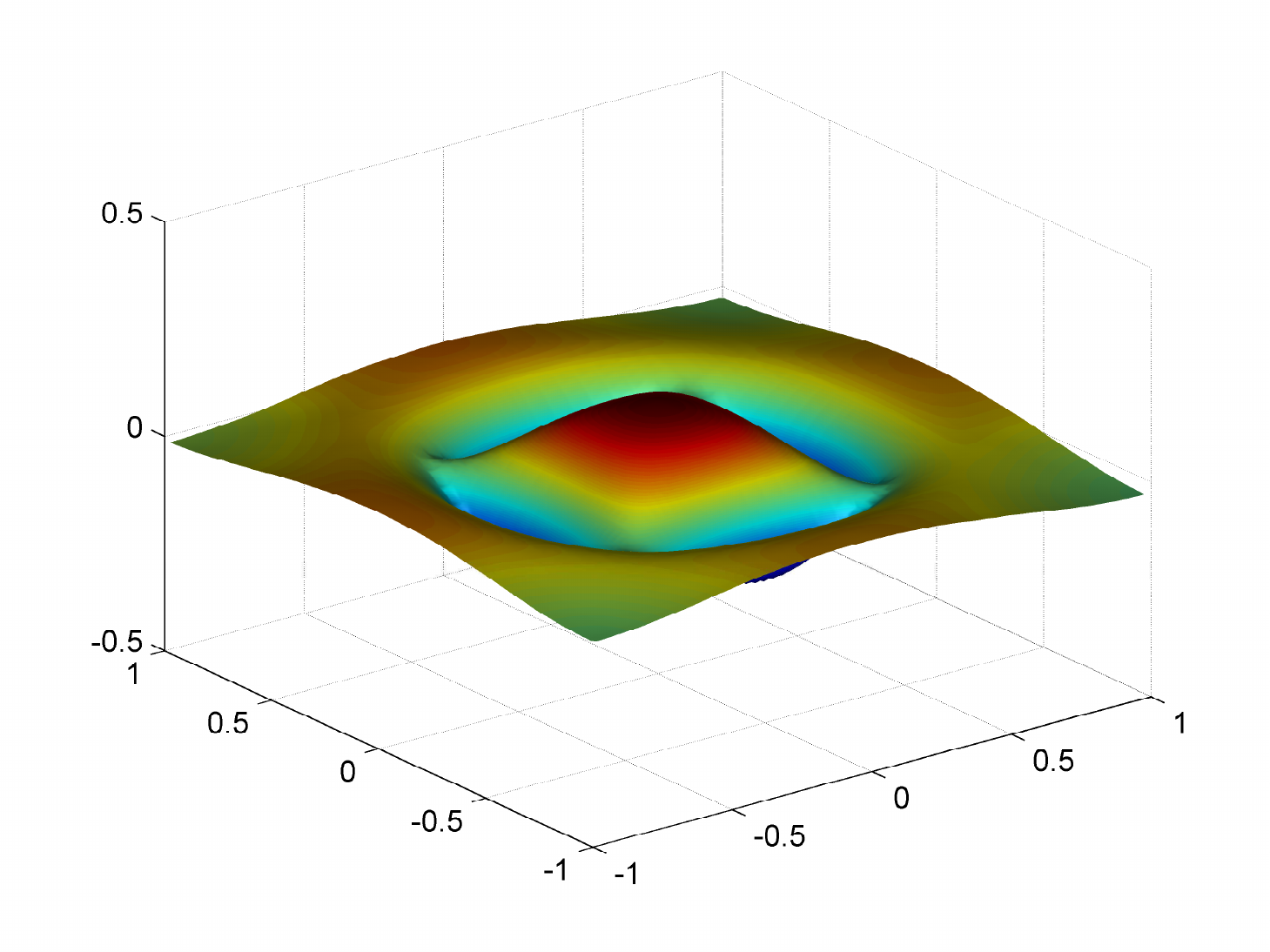} }
\subfigure[Difference by TGSS-DBTS]{
\label{figure 4-5}
\includegraphics[scale=0.34]{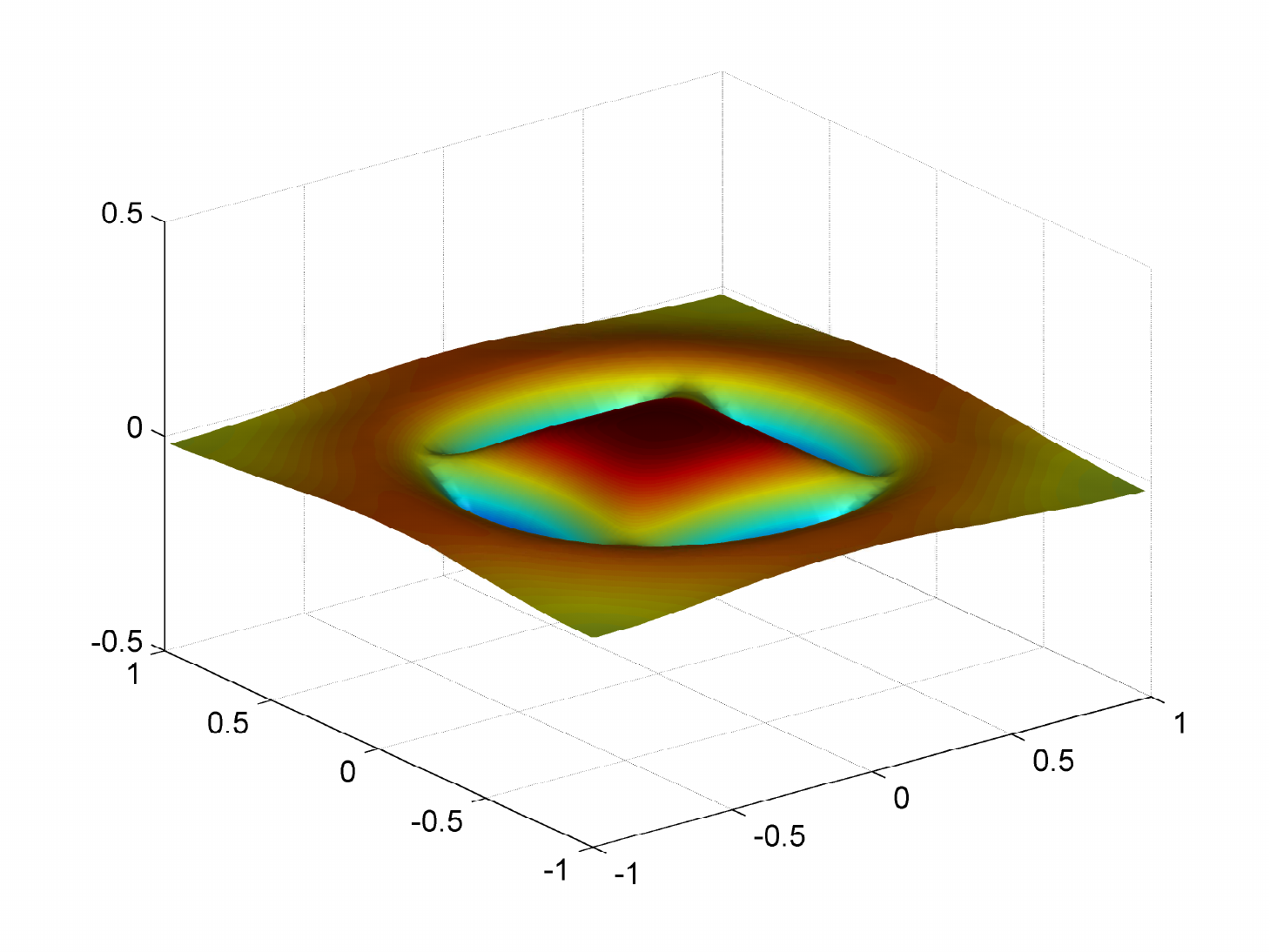} }
\subfigure[Difference by TGSS-Nes]{
\label{figure 4-6}
\includegraphics[scale=0.34]{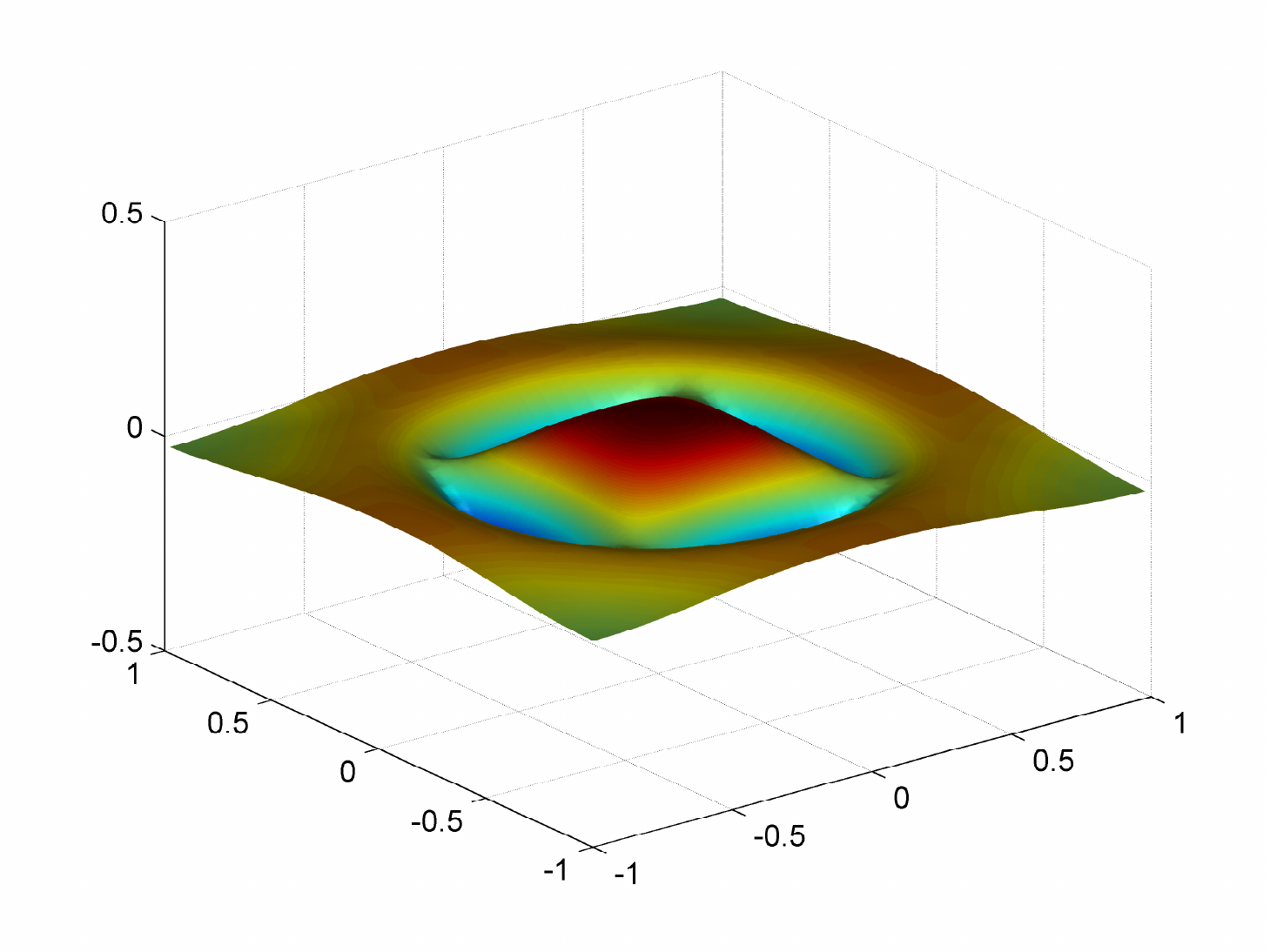} }
}
\caption{Difference $c^\dag-c_k$ by each method with $2\%$ noise level at 50th iteration.
}
\label{Noise-2D difference}
\end{figure}

For the case of noise data, we exhibit the RE curves in the first 100 iterations see Figure \ref{figure 2-3} and the reconstruction difference results at the 50th iteration for the noisy data ($\delta=2\%$) case in Figure \ref{Noise-2D difference} of different methods. We can intuitively see that at $\delta=2\%$, RE's value in TGSS method declines faster which means that it converges to the true solution faster. The more obvious effect is shown in Figure \ref{Noise-2D difference}. By comparing the Difference between exact and reconstructed solutions, we can observe that TGSS certainly reconstruct the coefficient of (\ref{equation 5-1}) better than TPG. Without losing generality, more detailed comparisons for noisy data are summarized in Table \ref{table 1-2}.

\begin{table}
\footnotesize
\centering
\caption{
Comparisons between Land, TPG, SESOP and TGSS methods for 2-D inverse potential problem with noisy data.
}
\begin{threeparttable}
\begin{tabular}{l l l l l l l l l l l l l}
\hline
$\delta$& Methods& $k_*$& Time (s)&  RE & Rate($k_*$) & & Rate(t)&\\
\hline
\multirow{1}{*}{2$\%$} & $\textrm{\textrm{Land}}$& 2451 & 199.289& $ {4.46 \times 10^{-2}}$ & 100$\%$ & & 100$\%$ & \\
&$\textrm{\textrm{TPG-DBTS}}$& 257 & 21.003 & ${4.48\times 10^{-2}}$ & 10.49$\%$ & & 10.54$\%$ & \\
&$\textrm{\textrm{TPG-Nes}}$& 205 & 16.679 & ${4.47 \times 10^{-2}}$ & 8.36$\%$ & & 8.37$\%$ & \\
&$\textrm{\textrm{SESOP}}$& 64 & 15.091 & ${4.14 \times 10^{-2}}$ & 2.61$\%$ & & 7.57$\%$ & \\
&$\textrm{\textrm{TGSS-DBTS}}$& 40 & 7.397 & ${3.88 \times 10^{-2}}$ & 1.63$\%$ & & 3.71$\%$ & \\
&$\textrm{\textrm{TGSS-Nes}}$& 40 & 6.138 & ${3.82 \times 10^{-2}}$ & 1.63$\%$ & & 3.08$\%$ & \\
\hline
\multirow{1}{*}{1$\%$} & $\textrm{\textrm{Land}}$& 4607 & 377.001 & ${2.82 \times 10^{-2}}$ & 100$\%$ & & 100$\%$ & \\
&$\textrm{\textrm{TPG-DBTS}}$& 376 & 30.845 & ${2.84 \times 10^{-2}}$ & 8.16$\%$ & & 8.18$\%$ & \\
&$\textrm{\textrm{TPG-Nes}}$& 274 & 22.400 & ${2.87 \times 10^{-2}}$ & 5.95$\%$ & & 7.42$\%$ & \\
&$\textrm{\textrm{SESOP}}$& 103 & 24.243 & ${2.61\times 10^{-2}}$ & 2.24$\%$ & & 6.43$\%$ & \\
&$\textrm{\textrm{TGSS-DBTS}}$& 68 & 12.106 & ${2.76 \times 10^{-2}}$ & 1.48$\%$ & & 3.21$\%$ & \\
&$\textrm{\textrm{TGSS-Nes}}$& 63 & 9.22 & ${2.86 \times 10^{-2}}$ & 1.37$\%$ & & 2.45$\%$ & \\
\hline
\multirow{1}{*}{0.5$\%$} & $\textrm{\textrm{Land}}$& 7293 & 595.095 & ${2.08 \times 10^{-2}}$ & 100$\%$ &  & 100$\%$ & \\
&$\textrm{\textrm{TPG-DBTS}}$& 469 & 38.010 & ${2.12 \times 10^{-2}}$ & 6.43$\%$ & & 6.39$\%$ & \\
&$\textrm{\textrm{TPG-Nes}}$& 337 & 27.969 & ${2.15 \times 10^{-2}}$ & 4.62$\%$ & & 4.70$\%$ & \\
&$\textrm{\textrm{SESOP}}$& 164 & 39.141 & ${2.01 \times 10^{-2}}$ & 2.25$\%$ & & 6.58$\%$ & \\
&$\textrm{\textrm{TGSS-DBTS}}$& 98 & 17.235 & ${2.05 \times 10^{-2}}$ & 1.34$\%$ & & 2.90$\%$ & \\
&$\textrm{\textrm{TGSS-Nes}}$& 97 & 14.534 & ${2.06\times 10^{-2}}$ & 1.33$\%$ & & 2.44$\%$ & \\
\hline
\end{tabular}
\end{threeparttable}
\label{table 1-2}
\end{table}

In order to visibly illustrate the convergence behavior of each method in 2-D inverse potential problem we created a table with detailed data. It can be seen from Table \ref{table 1-2} that at each noise level, TGSS takes about four to six times overall computational time and the iteration numbers less than TPG. That is to say, the convergence rate of the proposed TGSS is significantly accelerated compared to that of TPG. Under the same discrepancy principle, we can get similar conclusions to 1-D case. As can be expected, TGSS shows the favorable robustness.

To sum up, from the above analysis, we conclude that the TGSS method achieves more satisfactory acceleration effects than TPG method in respects accordingly.


\section{Conclusions}\label{section-6}
In this paper we introduced an accelerated two-point gradient method i.e. TGSS, as a paradigm of new iterative regularization method. TGSS can be regarded as a hybrid between two-point gradient and sequential subspace optimization methods for solving nonlinear ill-posed problem $F(x)=y$. As the Landweber method, we first find intermediate variables $z_k^\delta$ by combination parameters $\lambda_k^\delta$ according to Nesterov acceleration scheme or DBTS method, then substitute the original search direction with the one respect to $z_k^\delta$. As the TPG method, it uses information from previous iteration, which can improve the search direction to a finite number. For this purpose, we defined stripes for the exact and the noisy data case, considering the nonlinearity of the operator $F$ by utilize a tangential cone condition. The new iteration is updated by projecting the $z_k^\delta$ onto the intersections of stripes with respect to the search directions.

In terms of main theoretical results, we have presented convergence property of the TGSS method, then we get the regularization theory. Meanwhile, we have the corresponding consistency theory, which is to select $\lambda_k^\delta$ by DBTS method. In section \ref{section-5}, we presented numerical results for the TGSS method. The numerical experiments validate that the TGSS has excellent acceleration effects compared with Land, TPG and SESOP methods for producing satisfying approximations.

\section*{Acknowledgments}
This work was supported by the National Natural Science Foundation of China under Grant No.11871180.

\section*{References}

\end{document}